\RequirePackage{fix-cm}
\documentclass{svjour3}                     
\smartqed  
\usepackage[dvipdfm]{graphicx} 
\usepackage{bmpsize}
\journalname{Journal of Scientific Computing}

\usepackage{times}
\usepackage{graphics,epsfig,color}
\usepackage{amsmath,amssymb}
\usepackage{mathrsfs,amsfonts}

\usepackage{tikz}
\usetikzlibrary{calc,decorations.pathmorphing,patterns}
\usepackage{siunitx}
\usepackage{url}
\usepackage[debug,
colorlinks=true,
linkcolor=blue,
filecolor=green,
citecolor=blue,
pdfpagemode=UseNone]
{hyperref}

\newtheorem{alg}[theorem]{Algorithm}

\newtheorem{prop}[theorem]{Proposition}
\newtheorem{lem}[theorem]{Lemma}
\newtheorem{rem}[theorem]{Remark}
\DeclareMathAlphabet{\itbf}{OML}{cmm}{b}{it}
 \DeclareMathAlphabet\mathbfcal{OMS}{cmsy}{b}{n}

\renewcommand{\hat}{\widehat}
\renewcommand{\tilde}{\widetilde}
\def\RR{\mathbb{R}}
\def\CC{\mathbb{C}}

\def\vx{{\vec{\itbf x}}}

\def\bU{{\itbf U}}

\def\bGa{{\boldsymbol{\Gamma}}}

\def\bv{{\itbf v}}

\def\bnu{{\boldsymbol{\nu}}}

\def\bY{{\itbf Y}}

\def\bV{{\itbf V}}

\def\bG{{\itbf G}}
\def\bT{\boldsymbol{\mathbb{T}}}

\def\bGa{{\boldsymbol{\Gamma}}}
\def\bLa{{\boldsymbol{\Lambda}}}

\def\cL{{\mathcal L}}
\def\cW{{\mathcal W}}
\def\mW{\mathbb{W}}

\def\cS{{\mathcal S}}

\def\cQ{{\mathcal Q}}
\def\cR{{\mathcal R}}

\def\cZ{\mathcal{Z}}

\def\RM{{\scalebox{0.5}[0.4]{ROM}}}

\def\ep{\varepsilon}

\def\lb{\left <}
\def\rb{\right >}

\def\la{\lambda}

\begin{document}

\title{A reduced order model approach to inverse scattering in lossy layered media
}


\author{Liliana Borcea         \and
        Vladimir Druskin \and
         J\"{o}rn Zimmerling
}


\institute{Liliana Borcea   \at
             Department of Mathematics,\\
             University of Michigan,\\
  		  Ann Arbor, MI 48109-1043 \\
              \email{(borcea@umich.edu)}           
           \and
           Vladimir Druskin  \at
             Department of Mathematics,\\
             Worcester Polytechnic Institute,\\
  		 Worcester, MA 01609\\
              \email{(vdruskin@wpi.edu)}  
              \and
              J\"{o}rn Zimmerling
               \at
               {\it corresponding author}\\
             Department of Mathematics,\\
             University of Michigan,\\
  		  Ann Arbor, MI 48109-1043 \\
              \email{(jzimmerl@umich.edu)}  
}

\date{Received: date / Accepted: date}

\maketitle

\begin{abstract}
We introduce a reduced order model (ROM) methodology for  inverse electromagnetic wave scattering in layered lossy media, using data gathered by an antenna which generates a  probing wave 
and measures the time resolved reflected wave. We recast  the wave propagation problem 
as a passive infinite-dimensional dynamical system, whose transfer function is expressed in terms of the measurements at the antenna. The ROM is a low-dimensional  dynamical system that  approximates this transfer function.  While there are many possible ROM realizations, we are interested in one that preserves passivity   and in addition is: 
(1)  data driven (i.e., is constructed only from the measurements) and (2) it consists of a matrix with special sparse algebraic structure, whose entries contain spatially localized information about the unknown dielectric permittivity and electrical conductivity of the layered medium. Localized means in the intervals of a special  finite difference grid. The main result of the paper is to show with analysis and numerical simulations that these unknowns can be extracted efficiently from the ROM.

\end{abstract}

\section{Introduction}
\label{sect:intro}
We study a reduced order model (ROM) approach to inverse scattering  for Maxwell's equations in a lossy layered medium. We begin in section 
\ref{sect:1} with the formulation of the problem and show that it reduces to an inverse problem for a  stable and passive  dynamical system \cite{sorensen2005passivity,willems2007dissipative}, parametrized in port-Hamiltonian (pH) form \cite{beattie2019robust,jacob2012linear}. pH systems  arise from  energy based modeling of physical problems \cite{van2004port} and   are an active research area in reduced order modeling, where the focus is mainly on constructing a low-dimensional dynamical system, a ROM,  that approximates the transfer function of the original problem, while preserving the  pH structure 
 \cite{gugercin2012structure,benner2020identification}.  Because we are considering an inverse problem,  we develop  a data driven ROM which can be computed using simple numerical linear algebra tools that were developed originally for circuit synthesis \cite{Morgan2019ReflectionlessFT}.  The ROM looks like  a ladder RCL network model \cite{gugercin2012structure}. 
 It is described by a matrix of special sparsity pattern, corresponding to a finite difference scheme, with entries that  depend locally (on some grid) on the unknown dielectric permittivity and electrical conductivity. A similar ROM 
has been used for solving inverse scattering problems in lossless media in   \cite {druskin2016direct,borcea2019reduced}. However, for lossy media the algebraic sparsity constraint 
on the ROM does not permit obtaining a pH realization, meaning that some of the network resistors may be negative, in spite of the ROM being passive.    The main result of the paper, as we outline in section \ref{sect:2} of the introduction, is that such  ROMs can still   be used for an efficient solution of the inverse scattering problem or equivalently,  the parameter identification of the pH dynamical system. 

Inverse scattering  in layered lossy  media, although with different types of measurements,  has been considered in many other studies e.g. \cite{yagle1989one,bruckstein1985differential,jaulent1982inverse}. It can also be formulated as a  quadratic inverse spectral problem analyzed for example in \cite{buterin2012inverse,freiling2001inverse}.  
These results are specialized to one dimensional media. The ROM based inversion methodology introduced in this paper has the advantage that it can be extended, in principle, to inverse scattering in multi-dimensional media.

\subsection{The inverse scattering  problem}
\label{sect:1}
The electromagnetic waves are modeled by the  electric field $\vec{\bf E}(\vx,t)$ and magnetic field $\vec{\bf H}(\vx,t)$, satisfying  Maxwell's equations
\begin{align}
\begin{pmatrix}
\ep(z) & 0 \\
0 & \mu \end{pmatrix} \partial_t \begin{pmatrix} \vec{\bf E}(\vx,t) \\ \vec{\bf H}(\vx,t) \end{pmatrix} =
\begin{pmatrix} -\sigma(z) & \nabla \times \\ - \nabla \times & 0 \end{pmatrix} \begin{pmatrix} \vec{\bf E}(\vx,t) \\ \vec{\bf H}(\vx,t) \end{pmatrix}  {-} \begin{pmatrix}
\vec{{\bf J}}(\vx,t) \\ 0 \end{pmatrix},  
\label{eq:F1}
\end{align}
at time $t \in \mathbb{R}$ and location $\vx =(x_1,x_2,z)$ in the half space domain $\Omega = \mathbb{R}^2 \times (-\infty,L]$, with perfect electrical conductor boundary condition 
\begin{equation}
\vec{\bf e}_z \times \vec{\bf E}(\vx,t)   =0, \qquad \vx = (x_1,x_2,L), ~~ t \in \mathbb{R}.
\label{eq:F2}
\end{equation}
Here we  introduced the coordinate system with axes along the orthonormal vectors $\vec{\bf e}_{x_1}, \vec{\bf e}_{x_2}$ and $\vec{\bf e}_z$,
with $\vec{\bf e}_z$ pointing in the direction of variation of the layered medium. The coordinate $z \in (-\infty,L]$  along  $\vec{\bf e}_z$  is called the range and the other two (cross-range) coordinates are $(x_1,x_2) \in \mathbb{R}^2$.   
\begin{figure}[t!]
\centering
\includegraphics[width=0.45\textwidth]{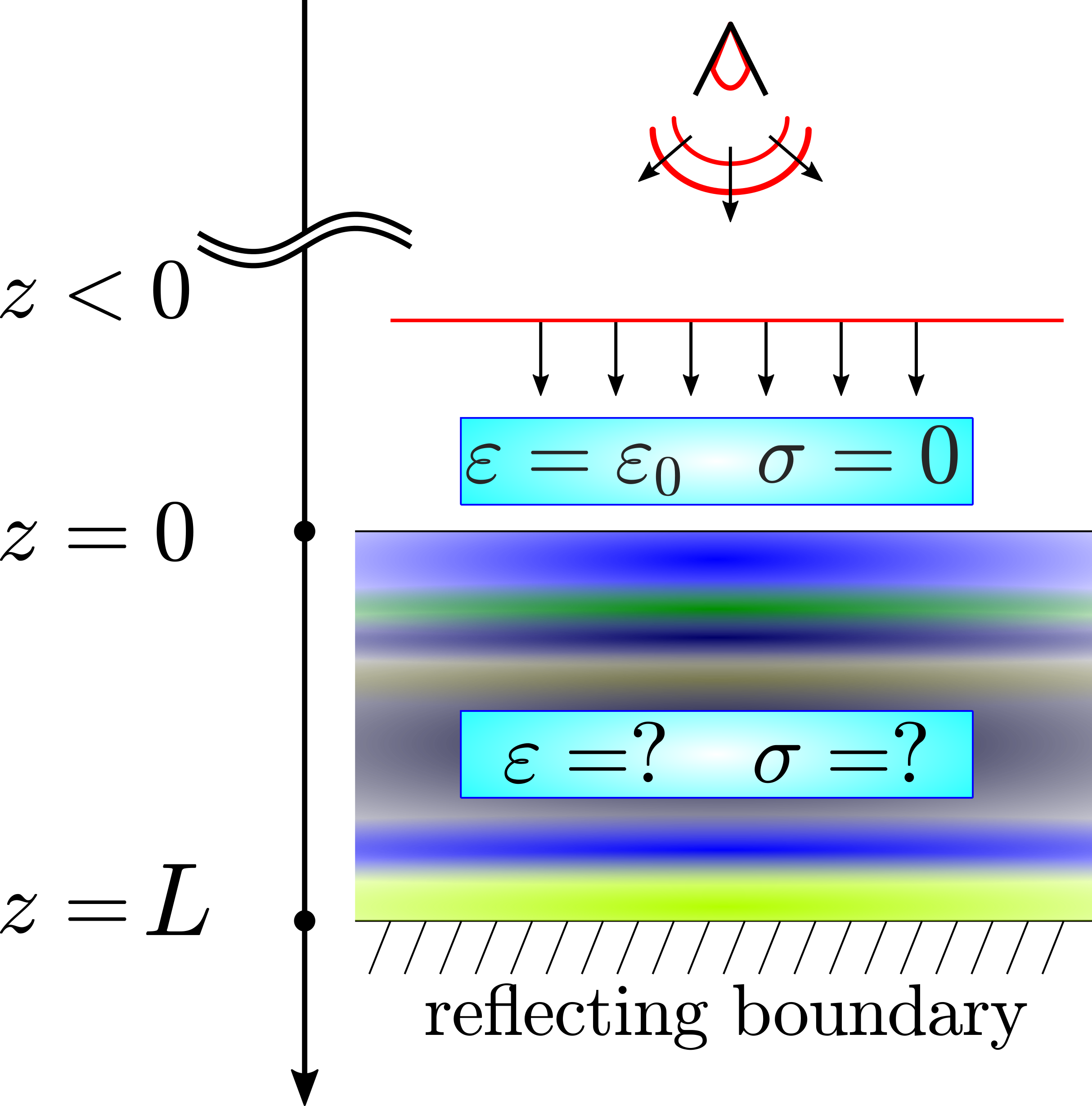}
\vspace{-0.0in}\caption{Illustration of the inverse scattering setup: A remote antenna  probes an unknown layered medium lying in the interval $z \in (0,L)$, with a downgoing wave that is approximately planar by the time 
it reaches the surface $z = 0$. The medium at $z < 0$ is known, lossless and homogeneous. The data for inverse scattering are the time resolved measurements of the upgoing wave. }
\label{fig:setup}
\end{figure}

Motivated by applications in radar imaging, we use the setup illustrated in Figure \ref{fig:setup}, where the medium is lossless and homogeneous at $z \in (-\infty, 0]$, with known dielectric permittivity $\ep(z) = \ep_0$, electrical conductivity $\sigma(z) = 0$ and magnetic permeablity $\mu$,
and it is variable and lossy at $z \in (0,L)$, with unknown positive $\ep(z)$ and  $\sigma(z)$. The magnetic permeability is assumed  constant throughout the domain. 

The data for the inverse problem are gathered by an antenna with directivity along the range direction, which emits waves and measures the 
scattered returns. The wave source is modeled in \eqref{eq:F1} by the  function $\vec{{\bf J}}(\vx,t)$ with short temporal support at time $t  \ge 0$  and with spatial support that is far away from the surface $z = 0$. 
We assume throughout that the wave is polarized so that the electric field points in the direction $-\vec{\bf e}_{x_2}$ and the magnetic field 
points in the direction $\vec {\bf e}_{x_1}$. Prior to the excitation there is no wave 
\begin{equation}
\vec{\bf E}(\vx,t) = \vec{\bf H}(\vx,t) = 0, \qquad \vx \in \Omega, ~~ t <0.
\label{eq:F3}
\end{equation}

Using causality and the finite wave speed,  we can calculate the incident wave impinging on the unknown layered medium  by solving Maxwell's equations with the source $\vec{{\bf J}}(\vx,t)$ in the whole space $\mathbb{R}^3$ filled with the homogeneous, lossless medium. Because the source is far away from the surface $z = 0$, 
for practical purposes we can approximate this wave  by a planar, down going wave with amplitude determined by the source. This wave penetrates the layered medium where it scatters, and the data for the inverse problem are the time resolved measurements of the reflected, upgoing wave captured at the antenna.  

The plane wave approximation allows us to reduce Maxwell's equations to a one dimensional inverse scattering problem  in the range interval $z \in (z_0,L)$, for some $z_0< 0 $ close to the interface $z=0$. Before writing this problem, let us transform to travel time coordinates 
\begin{equation}
T(z) = \int_0^z \sqrt{\mu \ep(z')} dz', 
\label{eq:TT}
\end{equation}
where we note that $1/\sqrt{\mu \ep(z)}$ is the wave speed. We also  
take  the Laplace transform with respect to $t$ and denote by $s \in \mathbb{C}$ the Laplace frequency. After these transformations, 
the component of the electric field along $- \vec{\bf e}_{x_2}$ denoted by 
$u(T,s)$,  aka the primary wave,  and  the component of the magnetic field along $\vec{\bf e}_{x_1}$ denoted by $\hat u(T,s)$, aka the dual wave,  satisfy the following system of coupled ordinary differential equations
\begin{align}
\partial_T \hat u(T,s) + \frac{\big[s + r(T)\big]}{\zeta(T)} u(T,s) &= 0, \label{eq:F4} \\
\partial_T u(T,s) + s \zeta(T) \hat u(T,s) &= 0, \label{eq:F5}
\end{align}
with homogeneous boundary condition at $T_L = T(L)$
\begin{equation}
u(T_L,s) = 0.
\label{eq:F6}
\end{equation}
The coefficients in these equations are the wave impedance 
\begin{equation}
\zeta(T) = \sqrt{\frac{\mu}{\ep(z(T))}},
\label{eq:F7}
\end{equation}
and the loss function
\begin{equation}
r(T) = \frac{\sigma(z(T))}{\ep(z(T))},
\label{eq:F8}
\end{equation}
the unknowns in the inverse problem. They satisfy $\zeta(T) = \zeta_0$ and $r(T) = 0$ for the negative travel time $T \in [T_0,0)$,
where $T_0 = T(z_0) < 0$.

We will restrict the problem to the interval $T \in [0,T_L]$.  To establish the boundary condition at $T = 0$, we use the forward (downgoing) and 
backward (upgoing) wave decomposition in the interval $(T_0,0)$ where the medium is homogeneous. This 
decomposition is 
\begin{equation}
\begin{pmatrix} 
u(T,s) \\ \hat u(T,s) \end{pmatrix} = 
\zeta_0^{\frac{1}{2}} \left[ a_f(T,s)e^{-sT} \begin{pmatrix} 1 \\ \zeta_0^{-1} \end{pmatrix} + 
a_b(T,s)e^{s T}\begin{pmatrix} 1 \\ -\zeta_0^{-1} \end{pmatrix} \right], \quad T \in (T_0,0),
\label{eq:F9}
\end{equation}
where $a_f(T,s)$ is the amplitude of the forward  going wave and $a_b(T,s)$ is the 
amplitude of the backward going wave. If we substitute \eqref{eq:F9} in equations \eqref{eq:F4}--\eqref{eq:F5}, we obtain that these amplitudes  satisfy
\begin{equation}
\partial_T a_f(T,s) = \partial_T a_b(T,s) = 0, \quad T \in (T_0,0).
\end{equation}
Therefore,  $a_f(0,s)$ is the amplitude of the computable incoming wave and $a_b(0,s)$ is the amplitude of the measured scattered wave, and \eqref{eq:F9}  gives
\begin{equation}
\frac{u(0,s)}{\hat u(0,s)} = \zeta_0 \left[ \frac{a_f(0,s) + a_b(0,s)}{a_f(0,s) - a_b(0,s)}\right].
\label{eq:F16}
\end{equation}
Our equations are linear, so we normalize henceforth the wave fields by setting 
\begin{equation}
\hat u(0,s) = 1.
\label{eq:F17}
\end{equation}

The  infinite dimensional dynamical system is described by the system of ordinary differential equations \eqref{eq:F4}--\eqref{eq:F5} for $T \in (0, T_L)$, with boundary conditions 
\eqref{eq:F6} and \eqref{eq:F17}. Its  transfer function $D(s) = u(0,s)$ is obtained from  \eqref{eq:F16} in terms of  the measurements at the antenna.
We assume we know $D(s)$ for  $s $ lying on some properly chosen curve in the complex plane, which depends on the signal emitted by the antenna.
 
\subsection{ROM approach to  system { parameter} identification}
\label{sect:2}
Let us rewrite the dynamical system in the following compact form 
\begin{equation}
\Big[ \cL + \cR(T) + s \cZ(T) \Big] \begin{pmatrix} u(T,s) \\ \hat u(T,s) \end{pmatrix} = \begin{pmatrix} \delta(T) \\ 0 
\end{pmatrix}, \qquad  T \in (0-,T_L),
\label{eq:F20}
\end{equation}
with homogeneous boundary conditions 
\begin{equation}
\begin{pmatrix} 1 &0 \\ 0 & 0 \end{pmatrix} \begin{pmatrix} u(T_L,s) \\ \hat u(T_L,s) \end{pmatrix} 
+ \begin{pmatrix} 0 &0 \\ 0 & 1 \end{pmatrix} \begin{pmatrix} u(0-,s) \\ \hat u(0-,s) \end{pmatrix} = {\bf 0}.
\label{eq:F20BC}
\end{equation}
Here $\cL$ is the skew-symmetric differential operator 
\begin{equation}
\cL = \begin{pmatrix} 0 & \partial_T \\ \partial_T & 0 \end{pmatrix}
\label{eq:F21}
\end{equation}
acting on vector valued  functions  satisfying the 
homogeneous boundary conditions \eqref{eq:F20BC}, and $\cR(T)$ and $\cZ(T)$ are the diagonal multiplication operators 
\begin{equation}
\cR(T)  = \begin{pmatrix} \frac{r(T)}{\zeta(T)} &0 \\ 0 & 0 \end{pmatrix}, \qquad 
\cZ(T)  = \begin{pmatrix} \frac{1}{\zeta(T)} &0 \\ 0 & \zeta(T) \end{pmatrix}.
\label{eq:F22}
\end{equation}
Note that $\cR(T)$ is positive semi-definite and $\cZ(T)$ is positive definite. 
\vspace{0.05in}
\begin{rem} In  \eqref{eq:F20} 
we have the slightly larger domain $z \in (0-,T_L)$, where $0-$ stands for a negative number, arbitrarily close to $0$. The boundary conditions \eqref{eq:F20BC} are homogeneous and we have a Dirac delta forcing which gives the jump condition\footnote{By $\hat u(0,s)$ we mean $\displaystyle \lim_{T \searrow 0} \hat u(T,s)$. }
\[\hat u(0,s) - \hat u(0-,s) = \hat u(0,s) = 1, \]
that  is consistent with the boundary condition \eqref{eq:F17}.  We use henceforth this formulation
with the Dirac delta forcing, instead of \eqref{eq:F4}--\eqref{eq:F6} and \eqref{eq:F17}, because it is easier to compare with the ROM dynamical system.
\end{rem}

\vspace{0.05in}

The inverse problem is:  
\emph{Determine the loss function $r(T)$ and the impedance function 
$\zeta(T)$  in the interval $T \in (0, T_L)$,  from the  transfer function }
\begin{equation}
D(s) = u(0,s) = \int_{0-}^{T_L} (\delta(T),0) \Big[ \cL + \cR(T) + s \cZ(T) \Big]^{-1} \begin{pmatrix} \delta(T) \\ 0 
\end{pmatrix} d T .
\label{eq:F23}
\end{equation}
We explain in  appendix \ref{ap:B.p}  that
\eqref{eq:F20}, with transfer function  \eqref{eq:F23} is a pH realization of a passive\footnote{Passive means that the dynamical system does not generate energy internally.} dynamical system. Therefore, we are solving a pH system identification  problem
\cite{benner2020identification}.

If we eliminate the dual wave $\hat u(T,s)$ in \eqref{eq:F20}-\eqref{eq:F20BC} we obtain an equivalent, 
second order formulation that is quadratic in the spectral parameter $s$. Therefore, the inverse problem can be 
recast as a quadratic inverse spectral problem that is uniquely solvable, as proved in \cite{buterin2012inverse}. 
An inversion method  is also given in \cite{buterin2012inverse} (with no numerical results), but it does not extend to multi-dimensional media. We introduce a novel ROM based approach to inversion, which can be generalized in principle to such media. The ROM is a low dimensional dynamical system that  looks like a finite difference scheme for the infinite dimensional dynamical system, with some caveats explained below.

If we had $\cR(T) = 0$, we could use the  data driven ROM considered in \cite{CPAM,druskin2016direct} 
\begin{equation}
\Big[ \boldsymbol{\cL} + s \boldsymbol{\cZ} \Big] \boldsymbol{\nu}(s) = \frac{{\bf e}_1}{\hat h_1},
\label{eq:F24}
\end{equation}
with $\boldsymbol{\cL},\boldsymbol{\cZ} \in \RR^{2n \times 2n},$ $\boldsymbol{\nu}(s) \in \CC^{2n}$ and 
${\bf e}_1$ the unit vector in $\RR^{2n}$ with the first entry equal to $1$ and all the other entries equal to $0$. 
The ROM is constructed so that its transfer function
\begin{equation}
D^{\RM}_n(s) = {\bf e}_1^T \Big[ \boldsymbol{\cL} + s \boldsymbol{\cZ} \Big]^{-1} \frac{{\bf e}_1}{\hat h_1},
\label{eq:F25}
\end{equation}
matches $2n$ measurements of $D(s)$. In \cite{CPAM} these measurements are the first $n$ poles and $n$ residues of $D(s)$ and the ROM is obtained via a layer-stripping algorithm based on the Lanczos recursion \cite{lanczos,chu2005inverse}. The ROM has special 
algebraic structure, with positive diagonal matrix $\boldsymbol{\cZ}$ and skew-symmetric bidiagonal matrix $\boldsymbol{\cL}$. Due to this structure, \eqref{eq:F24} can be interpreted as a finite difference discretization of \eqref{eq:F20} on a staggered grid with primary steps $h_j$  and dual steps $\hat h_j$, for $j = 1, \ldots, n$. This grid is obtained from the reference ROM constructed just as the data driven ROM, but for the reference medium with constant impedance. Its properties are described in \cite[Lemma 3.2]{CPAM}.
The entries in the data driven ROM are interpreted in \cite{CPAM} as local averages of the impedance on this grid, which leads to an easily computable  reconstruction $\zeta^{(n)}(T)$ indexed by $n$, that converges pointwise and in $L^1(0,T_L)$ to the true impedance $\zeta(T)$ in the limit $n \to \infty$ \cite[Theorem 6.1]{CPAM}.

The results in \cite{CPAM,druskin2016direct} were extended to  multi-dimensional media in \cite{DtB,borcea2019robust,borcea2019reduced}, where the ROM matrices $\boldsymbol{\cL}$ and 
$\boldsymbol{\cZ}$ have block bidiagonal and block diagonal structure.  It is not clear how to make a finite difference discretization analogy for such a ROM, because the blocks are full (non-sparse) matrices. Nevertheless, the block structure  captures the physics of wave propagation  in the range direction, 
and this has been exploited in \cite{borcea2019reduced} to devise a rapidly converging  optimization method 
for estimating the impedance function.

Our goal is to extend the ROM based inversion methodology developed in 
\cite{CPAM,druskin2016direct} to the identification of the coefficients of the pH system  \eqref{eq:F20}, where $\cR(T)$ is no longer zero, but positive semi-definite. As in \cite{CPAM}, we use the first 
$n$ poles and residues of the transfer function $D(s)$ to construct a $2n \times 2n$ matrix ROM\footnote{
These ``truncated spectral measure" measurements are used for convenience, 
but ROM's obtained from other matching conditions, such as $D^{\RM}_n(s_j) = D(s_j)$ for 
$2n$ properly chosen $(s_j)_{j=1}^{2n}$ can be used as well.}
\begin{equation}
\Big[ \boldsymbol{\cL} + \boldsymbol{\cR} + s \boldsymbol{\cZ} \Big] \boldsymbol{\nu}(s) = \frac{{\bf e}_1}{\hat h_1},
\label{eq:F26}
\end{equation}
with transfer function that has these poles and residues
\begin{equation}
D^{\RM}_n(s) = {\bf e}_1^T \Big[ \boldsymbol{\cL} + \boldsymbol{\cR}+ 
s \boldsymbol{\cZ} \Big]^{-1} \frac{{\bf e}_1}{\hat h_1}.
\label{eq:F27}
\end{equation}
Again, the ROM is obtained via direct layer stripping, implemented 
with the J-symmetric Lanczos algorithm \cite{Marshall1969synthesis,saad1982lanczos}.  
However,
there are two notable differences from the lossless case: 
(1) The Lanczos algorithm may break down, although it is unlikely to do so as explained in \cite{Joubert}. (2) The  tridiagonal ROM realization  may not be in pH form. That is to say, the diagonal 
$\boldsymbol{\cR}$ will have in general entries that model a non-existent ``magnetic loss" function $\hat r(T)$ that may take negative values. 
The latter difference makes the ROM based inversion more difficult than  in the lossless case,
as there is no direct analogy of  \eqref{eq:F26} as a finite difference scheme for the continuum problem \eqref{eq:F20}. 
If we tried to interpret it as a discretization for a problem with magnetic losses, we would arrive at an inverse problem for Telegrapher's (transmission line) equation  \cite{jaulent1982inverse,yagle1989one} that cannot be solved uniquely with our measurements. 
 The main difficulty addressed in this paper is how to embed the  ROM  \eqref{eq:F26} in the continuous pH
dynamical system \eqref{eq:F20}. We show with analysis and numerical simulations how to achieve this task for the case 
of losses with small enough amplitude variations. 

The paper is organized as follows: We begin in section \ref{sect:TF} with a discussion of the solvability of the inverse problem, the characterization of the transfer function 
and the truncated spectral measure data used for inversion. Then we give in section \ref{sect:ROM} the construction of the data driven ROM. In section \ref{sect:ROMinv} we introduce a simple, non-iterative  inversion algorithm,
and analyze it using first order perturbation analysis in the case of small variations of the loss function. 
We also show how to use optimization for dealing with losses with somewhat larger variation.
The analysis is complemented with numerical simulations. We end with a summary and a discussion of possible extensions to multi-dimensional media in
section \ref{sect:sum}.

\section{Solvability of the inverse problem and the  transfer function}
\label{sect:TF}
To describe the transfer function and to conclude that  it determines  
the impedance $\zeta(T)$ and the loss function $r(T)$ uniquely, we use the 
results in \cite{buterin2012inverse} on quadratic inverse spectral problems for 
Schr\"{o}dinger's equation with frequency dependent potential. To 
connect to these problems,  we assume henceforth\footnote{The regularity assumptions 
on $r(T)$ and $\zeta(T)$  can possibly be relaxed, but since we draw conclusions from 
\cite{buterin2012inverse}, we use the assumptions made in that study.} that $r(T)$ is absolutely continuous,
and that $\zeta(T)$ is smooth enough so that 
 \begin{equation}
q(T) = \zeta^{\frac{1}{2}}(T) \frac{d^2}{d T ^2} \zeta^{-\frac{1}{2}}(T) \in L^1([0,T_L]).
\label{eq:S6}
\end{equation} 
We also assume that the choice of the domain $(0-,T_L)$ is such that $\zeta(T)$ is constant in the small vicinity of $T  = 0$ and thus satisfies
\begin{equation}
\zeta(0) = \zeta_0, \qquad \frac{d\zeta}{d T } (0) = 0.
\label{eq:S1}
\end{equation}

The equation considered in \cite{buterin2012inverse}  is obtained from  \eqref{eq:F4}--\eqref{eq:F5} and  boundary conditions \eqref{eq:F6} and \eqref{eq:F17}, using the Liouville transformations 
\begin{equation}
w(T,s) = \sqrt{\frac{\zeta_0}{\zeta(T)}} u(T,s), \qquad \hat w(T,s) = \sqrt{\zeta_0 \zeta(T)} \hat u(T,s).
\label{eq:S5}
\end{equation}
Substituting these in \eqref{eq:F20}--\eqref{eq:F20BC} we obtain  the first order system 
\begin{align}
\left[ \cL + \cQ(T) + R(T) + s I\right] \begin{pmatrix} w(T,s) \\ \hat w(T,s) \end{pmatrix} = \begin{pmatrix} \zeta_0 \delta(T) \\ 0 
\end{pmatrix}, \qquad T \in (0-,T_L), 
\label{eq:SCF1}
\end{align}
with boundary conditions
\begin{equation}
\begin{pmatrix} 1 &0 \\ 0 & 0 \end{pmatrix} \begin{pmatrix} w(T_L,s) \\ \hat w(T_L,s) \end{pmatrix} 
+ \begin{pmatrix} 0 &0 \\ 0 & 1 \end{pmatrix} \begin{pmatrix} w(0-,s) \\ \hat w(0-,s) \end{pmatrix} = {\bf 0},
\label{eq:SCF2}
\end{equation}
where $I$ is the identity operator and $\cQ(T)$ and $R(T)$ are the multiplication operators
\begin{equation}
\cQ(T) = \begin{pmatrix} 0 & \frac{d}{d T} \ln \zeta^{-\frac{1}{2}}(T) \\- \frac{d}{d T} \ln \zeta^{-\frac{1}{2}}(T) & 0 
\end{pmatrix}, \qquad R(T) =  \begin{pmatrix} r(T) & 0\\ 0 & 0 
\end{pmatrix}.
\label{eq:SCF3}
\end{equation}
Eliminating the dual wave $\hat w(T,s)$ from these equations, we obtain an equivalent second order problem 
for the Schr\"{o}dinger equation  considered in \cite{buterin2012inverse}, 
\begin{align}
\big[\partial_T^2 - s^2 - s r(T) - q(T) \big] w(T,s) & = -s \zeta_0 \delta(T), \qquad T \in (0-,T_L), \label{eq:S2} \\
\partial_T w(0-,s) &= 0 , \quad 
w(T_L,s) = 0. \label{eq:S4}
\end{align}
The frequency dependent Schr\"{o}dinger   potential $s r(T) + q(T)$ is defined by the unknown loss function $r(T)$ and $q(T)$ given in \eqref{eq:S6}.  In an abuse of terminology, we shall refer to $q(T)$ as the potential. 
 
\subsection{The poles and zeroes of the transfer function}
\label{sect:TF1}
The result in \cite{buterin2012inverse}  is that the zeroes and poles of the so-called Weyl function $\cW(s)$
determine uniquely the potential $q(T)$ and the loss function $r(T)$. The impedance is then the solution of
the ordinary differential equation 
\begin{equation}
\frac{d^2}{d T ^2} \zeta^{-\frac{1}{2}}(T) = q(T) \zeta^{-\frac{1}{2}}(T), \qquad T \in (0,T_L), \label{eq:SQZ} 
\end{equation}
with initial conditions \eqref{eq:S1}.

We define the Weyl function in appendix \ref{ap:A} and show that it is related to the transfer function by 
\begin{equation}
D(s) = u(0,s) = w(0,s) =  -\frac{s\zeta_0 }{\cW(s)}.
\label{eq:S7}
\end{equation}
We also explain there  that if we introduce the quadratic (in $s$) Schr\"{o}dinger operator pencil 
\begin{equation}
\label{eq:A1}
\mathscr{L}_{q,r}(s) = \partial_T^2 - s^2 - s r(T) - q(T),
\end{equation}
then the poles of the transfer function, which are the zeroes of the Weyl function, are the eigenvalues 
of $\mathscr{L}_{q,r}(s)$ acting on the space $\cS_N$ of continuous functions satisfying 
\begin{equation} 
\phi(T) \in \cS_N ~~ {\rm if}~~ \partial_T \phi(0) = \phi(T_L) = 0.
\label{eq:S11}
\end{equation}
By an eigenvalue $\lambda \in \CC$ we mean that the null space of the  operator  is nontrivial \cite{markus2012introduction}. The zeroes of the transfer function are determined by the poles of the Weyl function, which are the eigenvalues of $\mathscr{L}_{q,r}(s)$ acting on the space 
$\cS_D$ of continuous functions satisfying 
\begin{equation} 
\phi(T) \in \cS_D ~~ {\rm if}~~ \phi(0) = \phi(T_L) = 0.
\label{eq:S10}
\end{equation}

The analysis in \cite{pronska2012spectral,pronska2013reconstruction} shows that the spectrum\footnote{The spectrum is defined as the set of $s \in \CC$ such that the operator is not boundedly invertible.} of the quadratic operator pencil \eqref{eq:A1} acting on either $\cS_N$  or $\cS_D$ coincides with the point spectrum
(i.e., all the spectral values are eigenvalues) and it is  countable. Moreover, definition \eqref{eq:A1} and the fact that  $r(T)$ and  $q(T)$ are real valued functions imply that these eigenvalues come in conjugate pairs. 

Consequently,  $D(s)$ is a meromorphic function, with 
poles  $\{\lambda_j, \overline{\lambda_j}, ~ j \ge 1 \}$ and  zeroes $\{\mu_j, \overline{\mu_j}, ~ j \ge 1 \}$, where the bar denotes the complex conjugate. The sets of poles and zeroes are disjoint, and they  determine explicitly the transfer function, as follows from the factorization of the Weyl function obtained in \cite{buterin2012inverse}. 
If we denote the mean loss by 
\begin{equation}
r_0 = \frac{1}{T_L} \int_0^{T_L} r(T) d T ,
\label{eq:VL2}
\end{equation}
then 
\begin{equation}
D(s) =  \frac{s \zeta_0 \Delta_D(s)}{\Delta_N(s)},
\label{eq:VL1}
\end{equation}
where the denominator in \eqref{eq:VL1} is given by \cite[Eq. (2.6)]{buterin2012inverse}
\begin{align}
\Delta_N(s) =& \frac{2}{r_0} \sinh \Big( \frac{T_L r_0}{2} \Big) \exp \Big\{
s T_L \Big[ \coth \Big( \frac{T_L r_0}{2}  \Big) - \frac{2}{T_L r_0 } \Big] \Big\} \nonumber \\
&\times \prod_{j=1}^\infty \frac{(s-\la_j)(s-\overline{\la_j})}{\big[\big(\frac{ \pi (j-1/2)}{T_L} \big)^2 + \frac{r_0^2}{4}\big]} \exp \left[ - \frac{s r_0}{\big(\frac{ \pi (j-1/2)}{T_L} \big)^2 + \frac{r_0^2}{4}} \right], 
\label{eq:VL3}
\end{align}
and the numerator is \cite[Eq. (2.14)]{buterin2012inverse}
\begin{align}
\Delta_D(s) =& \cosh \Big( \frac{T_L r_0}{2} \Big) \exp \Big[ s T_L \tanh \Big( \frac{T_L r_0}{2}  \Big) \Big] 
 \nonumber \\
&\times \prod_{j=1}^\infty \frac{(s-\mu_j)(s-\overline{\mu_j})}{\big[\big(\frac{ \pi j}{T_L} \big)^2 + \frac{r_0^2}{4}\big]} \exp \left[ - \frac{s r_0}{\big(\frac{ \pi j}{T_L} \big)^2 + \frac{r_0^2}{4}} \right].
\label{eq:VL4}
\end{align}
Moreover, we have the asymptotic expansion \cite[Eq. (2.5)]{buterin2012inverse} 
\begin{align}
\lambda_j &= \frac{i (j-1/2) \pi}{T_L} - \frac{r_0}{2} + \frac{i}{\pi(2j-1)} \int_0^{T_L} \Big[
q(T)-\frac{r^2(T)}{4} \Big] d T  + o\Big(\frac{1}{j} \Big)  \label{eq:asymptExp1}
\end{align}
which shows that $r_0$ can be determined from the real part of the poles in the limit $j \to \infty.$

\subsection{The truncated spectral measure transfer function of the ROM}
\label{sect:TF2}

We shall use the poles and residues representation of the transfer function $D(s)$, and assuming that the poles are simple, we can define the residues  $\{y_j, \overline{y_j}, ~j \ge 1\}$ by 
\begin{align}
y_j &= - \zeta_0 \la_j \Delta_D(\la_j)\lim_{s \to \la_j} \frac{(s-\lambda_j)}{\Delta_N(s)}.
\label{eq:resid}
\end{align}
The data driven ROM constructed in the next section has the ``truncated spectral measure" transfer function
\begin{equation}
D_n^{\RM}(s) = \sum_{j=1}^n \left[ \frac{y_j}{s - \lambda_j} + \frac{\overline{y_j}}{s - \overline{\lambda_j}} \right],
\label{eq:SL4}
\end{equation}
which shares the first $n$ poles $(\lambda_j)_{j=1}^n$ and residues $(y_j)_{j=1}^n$ with  $D(s)$.
We now  give a justification for this expression, in the case of a loss function $r(T)$ with small enough amplitude variations. 

\subsubsection{Constant loss function}
\label{sect:CLoss}
In the special case  $r(T) = r_0$,  the pencil \eqref{eq:A1} becomes 
\begin{equation}
\mathscr{L}_{q,r_0}(s) = \mathscr{L}_{q} - (s^2 + s r_0), \qquad \mathscr{L}_{q} = \partial_T^2 - q(T),
\label{eq:SL1}
\end{equation}
where $\mathscr{L}_{q}$ is a regular Sturm-Liouville operator. Since $\mathscr{L}_{q}$ acting on  $\cS_N$  has simple eigenvalues
$(-\theta_j^2)_{j=1}^\infty$, 
we conclude that the poles of the transfer function are 
simple in this case. In fact, we can write explicitly the series expansion of the transfer function, by 
expanding the solution of \eqref{eq:S2}--\eqref{eq:S4} in the orthonormal eigenbasis 
$\big(\varphi_j(T)\big)_{j = 1}^\infty$
of  $\mathscr{L}_{q}$.
We obtain that 
\begin{align}
D(s) &=  \int_{0-}^{T_L} \delta(T) \big[\mathscr{L}_{q,r_0}(s)\big]^{-1} \big[-s \zeta_0 \delta(T)\big] d T 
\nonumber \\
 &=\sum_{j=1}^\infty \frac{s \zeta_0\varphi_j^2(0)}{s^2 + s r_0 + \theta_j^2} = \sum_{j=0}^\infty 
\left[ \frac{y_j}{s - \lambda_j} + \frac{\overline{y_j}}{s - \overline{\lambda_j}} \right],
\label{eq:SL2}
\end{align}
where 
\begin{equation}
\lambda_j = -\frac{r_0}{2} + i \sqrt{\theta_j^2 - \frac{r_0^2}{4}}, \qquad y_j = \frac{\zeta_0 \varphi_j^2(0) \lambda_j}{\lambda_j - \overline{\lambda_j}}, \qquad 
j  \ge 1,
\label{eq:SL3}
\end{equation}

\subsubsection{Variable loss function}
\label{sect:VLTF}
Equations \eqref{eq:VL1}--\eqref{eq:VL4} define explicitly the transfer function in terms of the zeroes and poles
but it is difficult to obtain from it a series expansion like \eqref{eq:SL2}  for an arbitrary
non-negative loss function $r(T)$. In particular, the poles may no longer be simple. 
To avoid this complication, we assume henceforth that the variations of $r(T)$  about its mean 
$r_0 $ are not too large, so we can use the analytic perturbation theory of the eigenvalues of linear operators (see appendix \ref{ap:B}) to justify the series expression
\begin{align}
D(s) &=  \int_{0-}^{T_L} \delta(T) \big[\mathscr{L}_{q,r}(s)\big]^{-1} \big[-s \zeta_0 \delta(T)\big] d T 
= \sum_{j=0}^\infty 
\left[ \frac{y_j}{s - \lambda_j} + \frac{\overline{y_j}}{s - \overline{\lambda_j}} \right].
\label{eq:TRANSF}
\end{align}

We refer to appendix \ref{ap:B} for a proof that $D(s)$ satisfies the passivity conditions \cite{sorensen2005passivity,beattie2019robust,willems2007dissipative}.
The transfer function \eqref{eq:SL4} of the ROM approximates $D(s)$ by truncating this series after the 
$n^{\rm th}$ term. The error of the approximation at $s \notin \{\la_j, \overline \la_j, ~ j \ge 1\}$ is $O(1/n)$, as follows from equation \eqref{eq:resid} and the asymptotic expansion \eqref{eq:asymptExp1}.

\vspace{0.05in}
\begin{rem}
We will see in the next section that there exists a unique ROM of the form \eqref{eq:F26}, with transfer function 
\eqref{eq:SL4}. However, due to the imposed tridiagonal algebraic structure, the matrix $\boldsymbol{\cR}$ is not guaranteed
to be positive semi-definite, which means that the ROM does not preserve the pH structure of the dynamical system 
\eqref{eq:F20}. Nevertheless, the ROM preserves stability and all the numerical evidence is that it preserves passivity too (see section \ref{sect:ROMinv3}). The only reason we choose to work with the truncated spectral measure transfer function is because
we wish to make use of the analysis in \cite{CPAM}. However, there are other and likely better choices of data interpolation for constructing the ROM. A particularly interesting one is to interpolate the first $n$ poles and $n$ zeroes of the transfer function,  in which case the ROM is guaranteed to be passive \cite[Theorem 2.1]{sorensen2005passivity}. 
\end{rem}

\section{Construction and properties of the ROM}
\label{sect:ROM}
The discussion in section \ref{sect:TF} shows that we can formulate the inverse problem as:  
\emph{Given the poles $(\lambda_j)_{j=1}^n$ and residues $(y_j)_{j=1}^n$
of $D(s)$, obtain a ROM of the form \eqref{eq:F26}, with transfer function \eqref{eq:SL4}, and use it to compute   estimates of the impedance $\zeta(T)$ and the loss $r(T)$.} 
In this section we explain how to construct the ROM.

\subsection{The ROM as a finite difference scheme}
\label{sect:ROM.1}
We will obtain a ROM that corresponds to the  algebraic system
\begin{align}
\frac{\hat u_j(s) - \hat u_{j-1}(s)}{\hat h_j} + \frac{(s + \mathfrak{r}_j)}{\zeta_j} u_j(s) &= 0, \qquad j = 1, \ldots, n, \label{eq:R1}\\ 
\frac{u_{j+1}(s)-u_j(s)}{h_j} + (s + \hat{\mathfrak{r}}_j) \hat \zeta_j \hat u_j(s) &= 0, \qquad j = 1, \ldots, n, \label{eq:R2}\\
\hat u_0(s) &= 1, \quad 
u_{n+1}(s) = 0,  \label{eq:R4} 
\end{align}
and has the given transfer function $u_1(s) = D^{\RM}_n(s)$.
This ROM can be viewed as a  finite difference scheme on a staggered grid with primary steps $(h_j)_{j=1}^n$ and dual steps $ (\hat h_j)_{j=1}^n$.  We leave the grid unspecified for now, but  it turns out that among all possible staggered grids, there is a special one that can be computed and is useful for devising a ROM based inversion method. See sections \ref{sect:CONST2} and \ref{sect:ROMinv}.

 In \eqref{eq:R1}--\eqref{eq:R4} the primary wave is approximated by  $u_j$, at the grid points $T_j = \sum_{p=1}^{j-1} h_p$, and $\hat u_j$ approximates the dual wave at the  points $\hat T_j = \sum_{p=1}^j \hat h_p$, where   $T_1 = \hat T_0 = 0$. 
The coefficients $\zeta_j$ and $\hat \zeta_j$ model the impedance  on 
the staggered grid, and $\mathfrak{r}_j$ and $\hat{\mathfrak{r}}_j$ model losses. 

If we organize the wave approximations   in 
$\bU(s) = (u_1,\hat u_1, u_2, \ldots, u_n, \hat u_n )^T \in \CC^{2n},$ we can rewrite \eqref{eq:R1}--\eqref{eq:R4} in the $2n \times 2n$ matrix form  \eqref{eq:F26} given in the introduction.
The  bidiagonal  matrix
\begin{equation}
\boldsymbol{\cL} = {\rm diag} \big[(\hat h_1^{-1}, h_1^{-1}, \ldots, \hat h_n^{-1}),1\big] - 
{\rm diag} \big[(h_1^{-1}, \hat h_2^{-1}, \ldots, h_n^{-1}),-1],
\end{equation}
where ${\rm diag}[(\ldots),1]$ denotes the superdiagonal and ${\rm diag}[(\ldots),-1]$
the subdiagonal, is the  discrete analogue of the operator $\cL$ defined in \eqref{eq:F21}, 
and the diagonal  matrices 
\begin{equation}
\boldsymbol{\cR} = {\rm diag} (\mathfrak{r}_1/\zeta_1, \hat{\mathfrak{r}}_1 \hat \zeta_1, \ldots, \mathfrak{r}_n/\zeta_n, \hat{\mathfrak{r}}_n \hat \zeta_n), 
\qquad \boldsymbol{\cZ} = {\rm diag}(1/\zeta_1, \hat \zeta_1, \ldots, 1/\zeta_n, \hat \zeta_n), 
\end{equation}
look like  discretizations of the multiplication operators $\cR(T)$ and $\cZ(T)$ defined in 
\eqref{eq:F22}. However, the analogy is not  right because in $\boldsymbol{\cR}$ there are artificial dual losses $(\hat{\mathfrak{r}}_j)_{j=1}^n$ that may have negative values. This is the main difficulty in using the ROM for the inverse problem and will be 
addressed in section \ref{sect:ROMinv}. We explain next how to compute the ROM.

\subsection{Data driven ROM}
\label{sect:ROM2.p}
There are more unknowns in the system \eqref{eq:R1}--\eqref{eq:R4} than we can determine from the given
$D^{\RM}_n(s)$, so let us define  the coefficients that we can determine: 
\begin{align}
\gamma_j = h_j \hat \zeta_j, \quad \hat \gamma_j = \frac{\hat h_j}{\zeta_j}, \quad \mathfrak{r}_j ~ ~\mbox{and} ~~ \hat{\mathfrak{r}}_j , 
\quad \mbox{for}~ j = 1, \ldots, n.
\label{eq:R5}
\end{align}
We show in Appendix \ref{ap:C} that 
\begin{equation}
D^{\RM}_n(s) =  {\bf e_1}^T ({\bf A} + s {\bf I}_{2n})^{-1}  
\frac{{\bf e}_1}{\hat \gamma_1} = \sum_{j=1}^n \left[ \frac{y_j}{s - \lambda_j} + \frac{\overline{y_j}}{s - \overline{\lambda_j}} \right],
\label{eq:R9}
\end{equation}
where ${\bf I}_{2n}$ is the $2n \times 2n$ identity matrix and ${\bf A}$ is the complex symmetric, tridiagonal matrix 
\begin{align}
{\bf A} &= {\rm diag}(\alpha_1, \alpha_2, \ldots, \alpha_{2n}) + 
{\rm diag}\big[(\beta_2, \ldots, \beta_{2n}),1] + {\rm diag}\big[(\beta_2, \ldots, \beta_{2n}),-1],
 \label{eq:defA}
\end{align}
with off-diagonal entries 
\begin{align}
\beta_{2j} & = \frac{1}{\sqrt{-\gamma_j \hat \gamma_j} }, \quad j = 1, \ldots, n, \quad 
\beta_{2j+1} = - \frac{1}{\sqrt{-\gamma_j \hat \gamma_{j+1}} }, \quad j = 1, \ldots, n-1, \label{eq:R15} 
\end{align}
and diagonal entries 
\begin{align}
\alpha_{2j-1} &= \mathfrak{r}_j, \quad \alpha_{2j} = \hat{\mathfrak{r}}_j, \quad j = 1, \ldots, n.  \label{eq:R15a}
\end{align} 
We explain next how to calculate ${\bf A}$ from the given $(\la_j, y_j)_{j=1}^n$. 

\subsubsection{The Lanczos algorithm}
\label{sect:ROM2}
Let us organize the eigenvalues of ${\bf A}$, which are the negative of the poles in \eqref{eq:R9}, in the 
matrix  
\begin{equation}
\bLa = -\begin{pmatrix} 
\rm{diag}(\la_1, \ldots, \la_n) & 0 \\
0 & \rm{diag}(\overline{\la_1}, \ldots, \overline{\la_n}) 
\end{pmatrix}.
\label{eq:defLa}
\end{equation}
The J-symmetric Lanczos algorithm \cite{Marshall1969synthesis,saad1982lanczos}
uses the diagonalization  ${\bf A}  = \bY^T \bLa \bY,$ where $\bY = (\bY_1, \ldots, \bY_{2n} ) \in \CC^{2n \times 2n}$ is complex orthogonal i.e., $\bY^{-1} = \bY^T$. Substituting this in
\eqref{eq:R9} we get 
\begin{equation}
\bY_1^T (\bLa + s {\bf I}_{2n} )^{-1} \bY_1 = \hat \gamma_1 \sum_{j=1}^n \left[ \frac{y_j}{s - \lambda_j} + \frac{\overline{y_j}}{s - \overline{\lambda_j}} \right],
\end{equation}
and therefore, the first column of $\bY$ is 
\begin{equation}
\label{eq:bY1}
\bY_1 = \hat \gamma^{\frac{1}{2}}_1 \big( \sqrt{y_1},  
\ldots, \sqrt{y_n}, \sqrt{\overline{y_1}}, \ldots, \sqrt{\overline{y_n}}\big)^T.
\end{equation}
Moreover, since $\bY_1^T \bY_1 = 1$,  we can determine $\hat \gamma_1$ in terms of the residues,
\begin{equation}
\hat \gamma_1^{-1} = \sum_{j=1}^n (y_j + \overline{y_j})  = 2{\sum_{j=1}^n \mbox{Re}(y_j)}.
\label{eq:R12}
\end{equation}

The calculation of ${\bf A}$ is carried out by the Lanczos recursion, which equates both sides of $\bY {\bf A} = \bLa \bY$ column by column, starting from the first:

\vspace{0.05in}
\begin{alg}[From the poles and residues to ${\bf A}$]
\label{alg:Lancz1}
\begin{itemize}
\itemsep 0.04in
\item[]
\item \textbf{Input:} $(y_j, \la_j)_{j=1}^n$
\item \textbf{Initialization:} 
\begin{itemize}
\itemsep 0.03in
\item[] Compute $\bY_1$ as in \eqref{eq:bY1}, with $\hat \gamma_1$ from  \eqref{eq:R12}.
\item[] $\bv_1 = \bLa \bY_1$ 
\item[] $\alpha_1 = \bv_1^T \bY_1$
\item[] $\bv_1 = \bv_1 - \alpha_1 \bY_1$
\end{itemize}
\item \textbf{Processing steps:}
\begin{itemize}
\itemsep 0.02in
\item[] For $j = 2, \ldots, 2n$ do
\item[] \hspace{0.1in} $ \beta_j = \sqrt{\bv_{j-1}^T \bv_{j-1}}$ 
\item[] \hspace{0.1in} Breakdown if $\beta_j = 0$. Else compute
\item[] \hspace{0.1in} $\bY_j = \bv_{j-1}/\beta_j$
\item[] \hspace{0.1in} $\bv_j = \bLa \bY_j$ 
\item[] \hspace{0.1in} $\alpha_j = \bv_j^T \bY_j$ 
\item[] \hspace{0.1in} $ \bv_j = \bv_j - \alpha_j \bY_j - \beta_j \bY_{j-1}$
\item[] End for
\end{itemize}
\item \textbf{Output:} coefficients $(\alpha_j)_{j=1}^{2n}$, $(\beta_j)_{j=2}^{2n}$
\end{itemize}
\end{alg}

\vspace{0.05in}
Once we computed ${\bf A}$, we can determine the parameters \eqref{eq:R5} by solving equations 
\eqref{eq:R15}--\eqref{eq:R15a} as follows:

\vspace{0.05in}
\begin{alg}[From ${\bf A}$ to the ROM]
\label{alg:Lancz2}
\vspace{0.03in}
\begin{itemize}
\itemsep 0.04in
\item[]
\item \textbf{Input:} $(\alpha_j)_{j=1}^{2n}$, $(\beta_j)_{j=2}^{2n}$ and $\hat \gamma_1$ calculated from \eqref{eq:R12}.
\item \textbf{Processing steps:}
\begin{itemize}
\itemsep 0.02in
\item[] For $j = 1:n$ do 
\item[] \hspace{0.1in} $\gamma_j = - \frac{1}{\hat \gamma_j \beta_{2j}^2}$
\item[] \hspace{0.1in} $\mathfrak{r}_j = \alpha_{2j-1}$ 
\item[]  \hspace{0.1in} $\hat{\mathfrak{r}}_j = \alpha_{2j} $ 
\item[]  \hspace{0.1in} If $j < n$ compute $\hat \gamma_{j+1} = - \frac{1}{\gamma_j\beta_{2j+1}^2}$
\item[] End for
\end{itemize}
\item \textbf{Output:} ROM coefficients $(\gamma_j, \hat \gamma_j, \mathfrak{r}_j, \hat{\mathfrak{r}}_j)_{j=1}^n$
\end{itemize}
\end{alg}

\vspace{0.05in}

\subsection{Properties of the ROM}
\label{sect:ROM3}
As long as the Lanczos algorithm does not break down, which is almost always the case according to
 \cite{Joubert}, we obtain a unique set of ROM coefficients $(\gamma_j, \hat \gamma_j, \mathfrak{r}_j, \hat{\mathfrak{r}}_j)_{j=1}^n$.
 Thus, the data driven ROM exists, but in order for it to have an interpretation as the finite difference scheme
 \eqref{eq:R1}--\eqref{eq:R4}, we would like $(\gamma_j, \hat \gamma_j)_{j=1}^n$ to be positive and 
 $(\mathfrak{r}_j, \hat{\mathfrak{r}}_j)_{j=1}^n$ to be real valued.

A simple induction argument gives that the vectors $\bv_j$ and $\bY_j$ generated by Algorithm \ref{alg:Lancz1} are of the form 
\[
\bv_j = (v_{1,j}, \ldots, v_{n,j}, \overline{v_{1,j}}, \ldots, \overline{v_{n,j}})^T, \quad 
\bY_j = (Y_{1,j}, \ldots, Y_{n,j},\overline{Y_{1,j}}, \ldots, \overline{Y_{n,j}})^T.
\]
This implies that $(\alpha_j)_{j=1}^n$ and $(\beta_j^2)_{j=1}^n$ are real valued 
\begin{align*}
\alpha_j &= \bY_j^T \bLa \bY_j = -\sum_{l=1}^n \big( \la_l Y^2_{l,j} + \overline{\la _lY_{l,j}^2} \big) 
= - 2 \mbox{Re} \Big(\sum_{l=1} \la_l Y^2_{l,j}\Big), \qquad j = 1, \ldots, n,
\\
\beta_j^2 &= \bv_{j-1}^T \bv_{j-1} = \sum_{l = 1}^n \big( v_{l,j-1}^2 + \overline{v_{l,j-1}^2} \big) = 
2 \mbox{Re} \Big( \sum_{l = 1}^n v_{l,j-1}^2 \Big), \qquad j = 2, \ldots, n.
\end{align*}
However,  to get positive $(\gamma_j, \hat \gamma_j)_{j=1}^n$ from Algorithm \ref{alg:Lancz2} 
we need $(\beta_j^2)_{j=2}^n$ to be not only real, but negative. 
We show next that this holds for the case  of a constant loss. 
Since the Lanczos recursion
defines a continuous mapping from $(y_j, \la_j)_{j=1}^n$ to $(\beta_j^2)_{j=1}^n$, and since the mapping 
$r(T) \mapsto  (y_j, \la_j)_{j=1}^n$ is continuous as well, the conclusion extends to 
a loss function with small enough variations about a constant value. 
\subsubsection{The case of a constant loss}
\label{sect:CONST1}
Recall the explicit expression  \eqref{eq:SL2} of the transfer function in the case 
$r(T) = r_0$,
with poles and residues defined in \eqref{eq:SL3}. The next proposition describes the data driven ROM 
obtained by Algorithms \ref{alg:Lancz1}--\ref{alg:Lancz2}:

\vspace{0.05in} \begin{prop}
\label{prop.1}
The coefficients in the data driven ROM computed from the set of first $n$ poles and residues \eqref{eq:SL3} 
satisfy $\gamma_j, \hat \gamma_j > 0$, $\mathfrak{r}_j = r_0$ and $\hat{\mathfrak{r}}_j = 0$, for $j = 1, \ldots, n$. 
\end{prop}

\vspace{0.05in} \emph{Proof:} 
The Lanczos procedure for computing the ROM gives a unique answer, so it suffices to show 
that the ROM with the coefficients stated in the proposition has the given transfer function. 

Consider the discrete scheme 
obtained from \eqref{eq:R1} multiplied by $\zeta_j$ and $\eqref{eq:R2}$ multiplied by $\hat \zeta_j^{-1}$, 
where $\mathfrak{r}_j = r_0$ and $\hat{\mathfrak{r}}_j = 0$, 
\begin{align}
\frac{\hat u_j(s)-\hat u_{j-1}(s)}{\hat \gamma_j} + (s + r_0) u_j(s) &= 0, \qquad j = 1, \ldots, n, \label{eq:P1} \\
\frac{u_{j+1}(s) - u_j(s)}{\gamma_j} + s \hat u_j(s) &=0, \qquad j = 1, \ldots, n, \label{eq:P2} 
\end{align}
with boundary conditions \eqref{eq:R4}. Eliminating the dual wave, we obtain 
the $n \times n$ system 
\begin{align}
\frac{1}{\hat \gamma_j} \Big[ \frac{u_{j+1}(s) - u_j(s)}{\gamma_j} - \frac{u_{j}(s) - u_{j-1}(s)}{\gamma_{j-1}}
\Big] - s(s+r_0) u_j &=0, \qquad j = 2, \ldots, n, \label{eq:P3} \\
\frac{1}{\hat \gamma_1} \Big[ \frac{u_{2}(s) - u_1(s)}{\gamma_1}\Big] - s(s+r_0) u_1 &= -\frac{s}{\hat \gamma_1},
\label{eq:P4} \\
u_{n+1}(s) &= 0, \label{eq:P5}
\end{align}
written in matrix form 
\begin{equation}
\Big[ \bG - s(s+r_0) {\bf I}_n \Big] \begin{pmatrix}u_1(s) \\ \vdots \\ u_n(s) \end{pmatrix} = - \frac{s{\bf e}_1}{\hat \gamma_1}.
\label{eq:P6}
\end{equation}
Straightforward calculation shows that if $\hat \bGa = {\rm diag} (\hat \gamma_1, \ldots, \hat \gamma_n)$, then 
$
\tilde \bG  = \hat \bGa^{\frac{1}{2}}  \bG  \hat \bGa^{-\frac{1}{2}} 
$
is the tridiagonal, symmetric matrix with entries 
\begin{align}
\tilde{G}_{jl} = &(1-\delta_{j,1}) \Big[
-\frac{1}{\hat \gamma_j} \Big(\frac{1}{\gamma_j} + \frac{1}{\gamma_{j-1}}\Big) \delta_{j,l}+ 
\frac{\delta_{j+1,l} }{\gamma_j \sqrt{\hat \gamma_j \hat \gamma_{j+1}}} + \frac{\delta_{j-1,l} }{\gamma_{j-1}
\sqrt{\hat \gamma_{j} \hat \gamma_{j-1}}}  \Big] \nonumber \\
&+ \delta_{j,1} \Big[- \frac{\delta_{l,1}}{\hat \gamma_1 \gamma_1}  + \frac{\delta_{l,2}}{\gamma_1\sqrt{ \hat \gamma_1\hat \gamma_2}} \Big], \label{eq:P7p} \qquad j,l = 1, \ldots, n,
\end{align}
so we multiply  \eqref{eq:P6} on the left by $-\hat \bGa^{\frac{1}{2}}$ to get 
\begin{equation}
\Big[  s(s+r_0) {\bf I}_n - \tilde \bG  \Big] \hat \bGa^{\frac{1}{2}} \begin{pmatrix}u_1(s) \\ \vdots \\ u_n(s) \end{pmatrix} =  \frac{s{\bf e}_1}{\sqrt{\hat \gamma_1}}.
\label{eq:P7}
\end{equation}

We seek the entries of $\tilde \bG$ that give 
\begin{align}
u_1(s) &= \frac{s}{\hat \gamma_1} {\bf e}_1^T \Big[s(s+r_0) {\bf I}_n - \tilde{\bG} \Big]^{-1} {\bf e_1} = \sum_{j=1}^n \Big[ \frac{y_j}{s - \la_j} + \frac{\overline{y_j}}{s - \overline{\la_j}} \Big] \label{eq:EQ1}
\end{align}
with $(\la_j)_{j=1}^n$ and $(y_j)_{j=1}^n$ defined in \eqref{eq:SL3}. 
Recall from \eqref{eq:SL2} that 
\begin{equation}
 \frac{y_j}{s - \la_j} + \frac{\overline{y_j}}{s - \overline{\la_j}}  = \frac{s \xi_j}{s(s + r_0) + \theta_j^2},
\label{eq:P8}
\end{equation}
where $-\theta_j^2 <0 $ are the eigenvalues of the operator $\mathscr{L}_q$ defined in \eqref{eq:SL1} 
and $\xi_j > 0$ are defined in terms of the eigenfunctions of this operator evaluated at $T = 0$.
Therefore, the symmetric tridiagonal matrix $\tilde{\bG}$ satisfies
\begin{equation}
{\bf e}_1^T \Big[s(s+r_0) {\bf I}_n - \tilde{\bG} \Big]^{-1} {\bf e_1} = 
\hat \gamma_1 \sum_{j=0}^n \frac{\xi_j}{s(s + r_0) + \theta_j^2}.
\label{eq:P9}
\end{equation}
It  is negative definite, with eigenvalues $(-\theta_j^2)_{j=1}^n$ and orthonormal eigenvectors gathered 
as columns in the orthogonal matrix $\boldsymbol{\Phi}$. The first row of $\boldsymbol{\Phi}$ equals $\sqrt{\hat \gamma_1} 
(\sqrt{\xi_1}, \ldots, \sqrt{\xi_n})$ and since $\boldsymbol{\Phi} \boldsymbol{\Phi}^T = {\bf I}_n$, we get that 
$
\hat \gamma_1^{-1} = \sum_{j=1}^n \xi_j.$ We can now use the Lanczos recursion \cite{lanczos,chu2005inverse} for Jacobi type matrices to compute $\tilde \bG$ and then  obtain the positive coefficients $\gamma_j$ and $\hat \gamma_j$ from equations 
\eqref{eq:P7p} as shown in \cite[Section 2.2]{CPAM}. $\quad \Box$

\vspace{0.07in}
\begin{rem}
\label{rem:2}
Proposition \ref{prop.1} and the results in sections \ref{sect:CONST2} and \ref{sect:ROMinv}
imply that when $r(T) = r_0$, the coefficients  $(\beta_{j}^2)_{j=2}^{2n}$ 
defined in \eqref{eq:R15} lie on the negative real axis, at distance of order $(\pi n)/(2T_L)$ from the origin. 
In the case of variable $r(T)$ we know that $(\beta_{j}^2)_{j=2}^{2n}$ are real valued and due to the continuity of the mappings  $r(T) \mapsto (\la_j,y_j)_{j=1}^n \mapsto (\beta_j^2)_{j=2}^{2n}$, they will remain negative for small enough variations $|r(T)-r_0|$.
\end{rem}

\vspace{0.07in}
\begin{rem} \label{rem:3}
Proposition \ref{prop.1} shows that the discrete quadratic inverse spectral problem with the truncated measure spectral data has 
an exact solution given by the tridiagonal matrix $\tilde{\bG}$ and the exact loss $r_0$. We also show in section 
\ref{sect:ROMinv1} that the entries in $\tilde{\bG}$ determine an approximation of the impedance that converges pointwise to the true one
in the limit $n \to \infty$. If the loss function is not constant, then the data driven ROM exists and is unique, but it has nonzero dual losses $(\hat{\mathfrak{r}}_j)_{j=1}^n$. Consequently, the algebraic first order system \eqref{eq:R1}--\eqref{eq:R4} does not have an equivalent discrete quadratic inverse spectral problem formulation.  That is to say, there is 
no tridiagonal matrix $\tilde{\bG}$ and diagonal matrix of primary losses such that 
\[\frac{s}{\hat \gamma_1} {\bf e}_1^T \Big[s^2 {\bf I}_n + s \, {\rm diag}(r_1, \ldots, r_n)- \tilde{\bG} \Big]^{-1} {\bf e_1}=  \sum_{j=1}^n \Big[ \frac{y_j}{s - \la_j} + \frac{\overline{y_j}}{s - \overline{\la_j}} \Big].
\]
One could try to find such matrices via nonlinear optimization, but then the data fit will not be accurate, and the optimization will likely be difficult to carry out because the objective function is not convex. In section \ref{sect:ROMinv3}
we propose a better optimization approach which uses the data driven ROM calculated with Algorithms \ref{alg:Lancz1}--\ref{alg:Lancz2} to estimate efficiently the impedance and loss functions. 
\end{rem}
%
%

\subsubsection{The spectrally matched grid}
\label{sect:CONST2}
For the case of a homogeneous lossless medium  with constant impedance $\zeta(T) = 1$, 
the data driven ROM was calculated explicitly in \cite[Appendix A]{CPAM}. It is described by 
the positive coefficients $\gamma_j = h_j$, $\hat \gamma_j = \hat h_j$
and the zero loss coefficients $\mathfrak{r}_j = \hat{\mathfrak{r}}_j = 0$, for $j = 1, \ldots, n$.  Moreover, according to \cite[Lemma 3.2]{CPAM}, 
\begin{equation}
\hat h_1 < h_1 < \hat h_2 < h_2 < \ldots < \hat h_n < h_n,
\label{eq:MONOT}
\end{equation}
and for sufficiently large $n$,
\begin{align}
\frac{h_j}{T_L} &= \frac{2 + O\big[ (n-j)^{-1} + j^{-2}\big]}{\pi \sqrt{n^2-j^2}}, \quad j = 1, \ldots, n-1, \quad 
\frac{h_n}{T_L} = \frac{\sqrt{2} + O(n^{-1})}{\sqrt{\pi n}}, \\
\frac{\hat h_j}{T_L} & = \frac{2 + O\big[ (n+1-j)^{-1} + j^{-2}\big]}{\pi \sqrt{n^2-(j-1/2)^2}}, \quad j = 1, \ldots, n.
\end{align}

\begin{figure}[t!]
\centering
\includegraphics[clip,trim=10mm 5mm 10mm 10mm,width=0.6\textwidth]{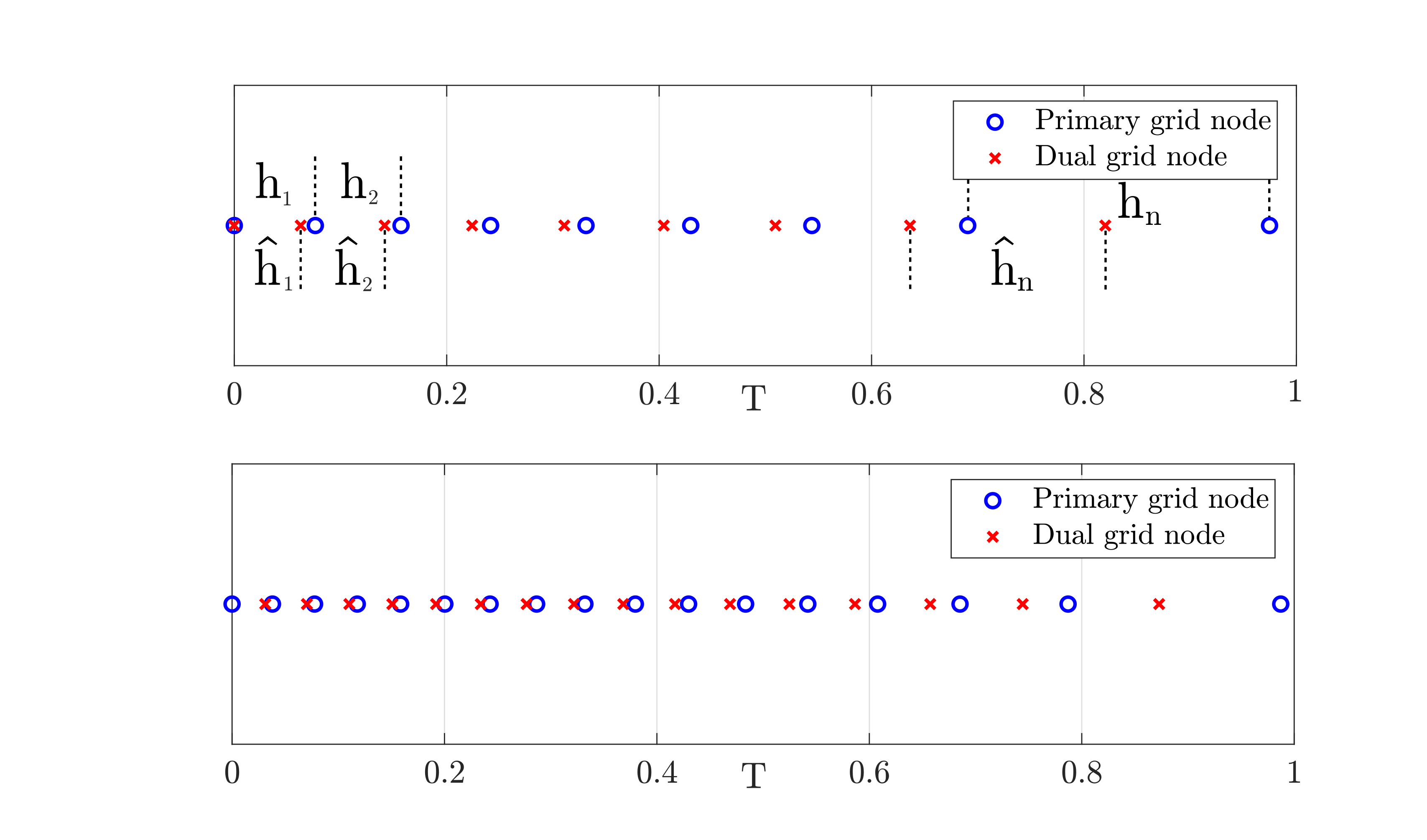}
\vspace{-0.0in}\caption{We display two spectrally matched grids calculated in the interval $T \in (0,T_L)$ 
with normalized $T_L = 1$. The top plot is for $n = 8$ and we indicate a few of the grid steps. The 
bottom plot is for $n = 16$.}
\label{fig:grid}
\end{figure}

Recalling equations \eqref{eq:R1}--\eqref{eq:R4}, we see that  $(h_j, \hat h_j)_{j=1}^n$ can be interpreted 
as grid steps for a finite difference scheme on the staggered grid with primary points 
\begin{align}
T_j &= \sum_{p=1}^{j-1} h_p, \qquad j = 2, \ldots, n+1, ~~ T_1 = 0,
\label{eq:prGP}
\end{align}
that are interlaced with the  dual points 
\begin{align}
\hat T_j &= \sum_{p=1}^j \hat h_j, \qquad j = 1, \ldots, n, ~~ \hat T_0 = 0.\label{eq:duGP}
\end{align}
See Figure \ref{fig:grid} for an illustration.  Since the finite difference scheme matches exactly the truncated spectral measure transfer function, we call the grid spectrally matched\footnote{In \cite{CPAM} the grid was called ``optimal", but spectrally matched is a more appropriate name.}. 
We will use this  grid in the scheme \eqref{eq:R1}--\eqref{eq:R4} and explain next why it
is the right choice. 

\section{ROM based inversion}
\label{sect:ROMinv}
We now show how to use the data driven ROM for solving the inverse problem. We begin 
with the easy case of a constant loss $r(T) = r_0$ in section \ref{sect:ROMinv1}. Then we consider in section \ref{sect:ROMinv2} a slightly varying $r(T)$ and propose a simple inversion algorithm that 
can be analyzed using first order spectral perturbation analysis. The main purpose of studying these two cases is to show explicitly that the ROM coefficients contain information about the unknown $r(T)$ and $\zeta(T)$ that is localized in the spectrally matched grid intervals. The conclusion extends to larger variations of $r(T)$, as shown with numerical simulations in section \ref{sect:ROMinv3}. 

\subsection{Inversion for a constant  loss function}
\label{sect:ROMinv1}

We described in Proposition \ref{prop.1} the data driven ROM for  the medium with constant loss $r(T) = r_0$. 
We saw that its coefficients satisfy $\gamma_j, \hat \gamma_j > 0$, $\mathfrak{r}_j = r_0$ and $\hat{\mathfrak{r}}_j = 0$, for $j = 1, \ldots, n$.  
Thus, we can determine  the constant loss exactly in this case. The next proposition shows that the impedance function 
is also easily estimated from the ROM:

\vspace{0.05in}
\begin{prop}
\label{prop.2}
Consider the data driven ROM computed with Algorithms \ref{alg:Lancz1}--\ref{alg:Lancz2} from the first $n$ poles and residues given in \eqref{eq:SL3}.
Define the coefficients $(\zeta_j,\hat \zeta_j)_{j=1}^n$ as in equation \eqref{eq:R5} using the spectrally 
matched grid steps $(h_j,\hat h_j)_{j=1}^n$. Let $\zeta^{(n)}(T)$ be some interpolation (e.g., piecewise constant or linear) of these 
coefficients on the spectrally matched grid  \eqref{eq:prGP}--\eqref{eq:duGP} i.e., 
\begin{equation}
\zeta^{(n)}(T_j) = \zeta_j, \quad \zeta^{(n)}(\hat T_j) = \hat \zeta_j, \qquad j = 1, \ldots, n.
\end{equation}
Then, 
$\zeta^{(n)}(T) \to \zeta(T)$ pointwise and in $L^1([0, T_L])$ as $n \to \infty$. 
\end{prop}

\vspace{0.05in} \emph{Proof:} This result follows  from the proof of Proposition \ref{prop.1},
where we explained that the coefficients $(\gamma_j,\hat \gamma_j)_{j=1}^n$ are exactly the same 
as those calculated in \cite[Section 2.2]{CPAM}. Then, we can cite directly  \cite[Theorem 6.2]{CPAM} which is in terms of the function $\zeta^{-1}(T)$. That theorem says that the piecewise 
constant interpolation of the values $\zeta_j^{-1} = \hat \gamma_j/\hat h_j$ and $\hat \zeta_j^{-1} = h_j/\gamma_j$
for $j = 1, \ldots, n$, converges pointwise and in $L^1([0, T_L])$ to $\zeta^{-1}(T)$, as $n \to \infty$.
$~\Box$

\subsection{Inversion for small variations of the loss function}
\label{sect:ROMinv2}
The ROM is constructed from $D^{\RM}_n(s)$ that has the first  $n$ poles and residues of $D(s)$. Recalling the expression
\eqref{eq:TRANSF} of $D(s)$ and \eqref{eq:SL4} of  $D^{\RM}_n(s)$, we have pointwise, for $s \in \CC \setminus \{\la_j , \overline{\la_j}, ~ j \ge 1\}$, that
\begin{equation}
\lim_{n \to \infty} 
D_n^{\RM}(s) = \lim_{n \to \infty} \sum_{j=1}^n \left[ \frac{y_j}{s - \lambda_j} + \frac{\overline{y_j}}{s - \overline{\lambda_j}} \right] 
= D(s).
\end{equation}
Thus, if $s$ is a zero of $D(s)$, we must have $D^{\RM}_n(s) \approx 0 $ for $n \gg 1$.
In this section we use a first order perturbation analysis of the poles and zeroes of $D(s)$ to explain how to 
devise a ROM based inversion algorithm on the spectrally matched grid defined in section \ref{sect:CONST2}.

To state the result, we need  the spectral decomposition 
of the linear operator 
\begin{equation}
\cL_\zeta = \begin{pmatrix} 0 & \zeta(T) \partial_T \\ \zeta^{-1}(T) \partial_T
\end{pmatrix},  
\label{eq:pencilP}
\end{equation}
acting on vector valued functions $\bnu(T)$ satisfying  either the boundary conditions 
\begin{equation}
\begin{pmatrix} 0 & 0 \\ 0 & 1 \end{pmatrix} \bnu(0) = \begin{pmatrix} 1 & 0 \\ 0 & 0 \end{pmatrix} \bnu(T_L) = {\bf 0},
\label{eq:BC1} 
\end{equation}
or 
\begin{equation}
\begin{pmatrix} 1 & 0 \\ 0 & 0 \end{pmatrix} \bnu(0) =  \begin{pmatrix} 1 & 0 \\ 0 & 0 \end{pmatrix} \bnu(T_L) = {\bf 0}.
\label{eq:BC2} 
\end{equation}

\begin{lem}
\label{lem.1}
The spectrum of the linear operator $\cL_\zeta$ with boundary conditions \eqref{eq:BC1}  consists of 
purely imaginary eigenvalues $\{\pm i \theta_j, ~j \ge 1\}$, where $\theta_j$ are the same as in section \ref{sect:CLoss}. The 
eigenfunctions are 
\begin{equation}
\boldsymbol{\Phi}_j^{\pm}(T) =  \frac{1}{\sqrt{2}}\begin{pmatrix} \phi_j(T) \\ \pm i \hat \phi_j(T) \end{pmatrix}, \qquad j \ge 1, 
\label{eq:eigPhi}
\end{equation}
where $\phi_j(T)$ and $\hat \phi_j(T)$ are real valued functions satisfying the orthogonality relations 
\begin{equation} 
\int_0^{T_L} \zeta^{-1}(T) \phi_j(T) \phi_l(T) d T = \int_0^{T_L} \zeta(T) \hat \phi_j(T) \hat \phi_l(T) d T = \delta_{jl}, 
\qquad j , l \ge 1.
\label{eq:ORT1}
\end{equation}
Similarly, the spectrum of the linear operator  $\cL_\zeta$ with boundary conditions \eqref{eq:BC2} consists of the purely  imaginary eigenvalues $\{\pm i \vartheta_j, ~j \ge 1\}$ and its eigenvectors are 
\begin{equation}
\boldsymbol{\Psi}_j^{\pm}(T) =  \frac{1}{\sqrt{2}}\begin{pmatrix} \psi_j(T) \\ \pm i \hat \psi_j(T) \end{pmatrix}, \qquad j \ge 1, 
\label{eq:eigPsi}
\end{equation}
with components defined by the real valued functions $\psi_j(T)$ and $\hat \psi_j(T)$ satisfying
\begin{equation} 
\int_0^{T_L} \zeta^{-1}(T) \psi_j(T) \psi_l(T) d T = \int_0^{T_L} \zeta(T) \hat \psi_j(T) \hat \psi_l(T) d T = \delta_{jl}, 
\qquad j , l \ge 1.
\label{eq:ORT2}
\end{equation}
\end{lem}

\vspace{0.05in} The proof of this lemma is a simple calculation given in Appendix \ref{ap:D}.  Now let us denote by $\zeta^{(n)}(T)$ the ROM based estimate of the impedance function,  defined by an  interpolation of the values $(\zeta_j,\hat \zeta_j)_{j=1}^n$ assigned to the  grid points $(T_j, \hat T_j)_{j=1}^n$ defined in \eqref{eq:prGP}--\eqref{eq:duGP}. Let also $\mathfrak{r}^{(n)}(T)$ be an interpolation of the ROM parameters $(\mathfrak{r}_j)_{j=1}^n$ assigned to the 
primary grid points $(T_j)_{j=1}^n$ and $\hat{\mathfrak{r}}^{(n)}(T)$ an interpolation of $(\hat{\mathfrak{r}}_j)_{j=1}^n$ assigned to the dual grid points $(\hat T_j)_{j=1}^n$. We choose the simplest interpolation: piecewise linear for the impedance and piecewise constant for the losses.  However, one can use other interpolations that result in smoother estimates, as assumed in section \ref{sect:TF}. 
The next proposition describes the relation between the estimates $\zeta^{(n)}(T)$, ${\mathfrak{r}}^{(n)}(T)$ and $\hat{\mathfrak{r}}^{(n)}(T)$ and the true impedance and loss functions.

\vspace{0.05in}
\begin{prop}
\label{prop.3}
Suppose that the loss function satisfies 
\begin{equation}
r(T) = r_0 + \alpha \rho(T), \qquad 
\sup_{T \in (0,T_L)}{|\rho(T)|}/{r_0} = O(1),  \qquad 0 < \alpha \ll 1.
\label{eq:modelr}
\end{equation}
Then, we have the following pointwise ROM based estimate of the impedance
\begin{equation}
\zeta^{(n)}(T) = \zeta(T)[1 + o(1) + O(\alpha^2)].
\label{eq:convZeta} 
\end{equation}
Moreover, the functions $\mathfrak{r}^{(n)}(T)$ and $\hat{\mathfrak{r}}^{(n)}(T)$ are of the form
\begin{equation} 
\mathfrak{r}^{(n)}(T) = r_0 + \alpha \rho^{(n)}(T)[ 1+ o(1) + O(\alpha)], \quad \hat{\mathfrak{r}}^{(n)}(T) = \alpha \hat \rho^{(n)}(T)[ 1+o(1) + O(\alpha)], 
\label{eq:convR}
\end{equation}
where the $O(\alpha)$ terms satisfy
\begin{align}
\hspace{-0.05in} \int_0^{T_L} \rho(T) \frac{\phi_j^2(T)}{\zeta_j(T)} d T &= \int_0^{T_L} \Big[\rho^{(n)}(T) \frac{\phi_j^2(T)}{\zeta_j(T)} + 
\hat \rho^{(n)}(T) \zeta(T) \hat \phi^2_j(T) \Big] dT, \label{eq:convR1} \\
\hspace{-0.05in} \int_0^{T_L} \rho(T) \frac{\psi_j^2(T)}{\zeta_j(T)} d T &= \int_0^{T_L} \Big[\rho^{(n)}(T) \frac{\psi_j^2(T)}{\zeta_j(T)} + 
\hat \rho^{(n)}(T) \zeta(T) \hat \psi^2_j(T) \Big] dT, \label{eq:convR2}
 \end{align}
for $j \ge 1$. Here $o(1)$ is in the limit $n \to \infty$  and we used the eigenfunctions  in Lemma \ref{lem.1}.
\end{prop}

\vspace{0.05in}
The proof of this proposition is given in Appendix \ref{ap:E}. The next inversion algorithm uses it to estimate the impedance and loss function.

\begin{alg}[Inversion on the spectrally matched grid]
\label{alg:invGrid}
\vspace{0.03in}
\begin{itemize}
\itemsep 0.04in
\item[]
\item \textbf{Input:} The ROM coefficients $(\gamma_j, \hat \gamma_j, \mathfrak{r}_j, \hat{\mathfrak{r}}_j)_{j=1}^n$ and the spectrally matched grid 
steps $(h_j,\hat h_j)_{j=1}^n$. 
\item \textbf{Estimate of the impedance function:}
\begin{itemize}
\itemsep 0.02in
\item Calculate $\zeta_j = \hat h_j/\hat \gamma_j$ and $\hat \zeta_j = \gamma_j/h_j$, for $j = 1, \ldots, n$. 
\item Define estimate $\zeta^{(n)}(T)$ as the piecewise linear interpolation of the values
\[\zeta^{(n)}(T_j) = \zeta_j, \quad \zeta^{(n)}(\hat T_j) = \hat \zeta_j, \qquad j = 1, \ldots, n.\]
\end{itemize}
\item \textbf{Estimate of the loss function:} 
\begin{itemize}
\itemsep 0.02in
\item Calculate the eigenfunctions $\boldsymbol{\Phi}_j^\pm(T)$ and $\boldsymbol{\Psi}_j^\pm(T)$ described in 
Lemma \ref{lem.1}, for the estimated impedance $\zeta^{(n)}(T)$.
\item Calculate 
\begin{align*}
\mathfrak{r}^{(n)}(T) &= \sum_{j=1}^n \mathfrak{r}_j 1_{[T_j,T_{j+1})}(T) + \mathfrak{r}_n 1_{[T_{n+1}, T_L]}(T),  \\
\hat{\mathfrak{r}}^{(n)}(T) &= \sum_{j=1}^n \hat{\mathfrak{r}}_j 1_{[\hat T_{j-1},\hat T_{j})}(T) + \hat{\mathfrak{r}}_n 1_{[\hat T_n, T_L]}(T),
\end{align*}
where $1_{[a,b)}(T)$ is the indicator function of the interval $[a,b)$, equal to $1$ for $T \in [a,b)$ and $0$ otherwise.
\item The estimate of the loss function is 
\[r^{(n)}(T) = \sum_{j=1}^n \big[ r_j1_{[T_j, \hat T_j)}(T) + \hat r_j 1_{[\hat T_j,T_{j+1})}(T)\big] + \hat r_n 1_{[T_{n+1}, T_L]}(T), \]
where $(r_j,\hat r_j)_{j=1}^n$ are obtained by solving the $2n \times 2n$ linear system 
\begin{align*}
\sum_{l=1}^n \Big[ r_l \int_{T_l}^{\hat T_l} \frac{\phi_j^2(T)}{\zeta_j(T)} d T + \hat r_l \int_{\hat T_l}^{T_{l+1}} \frac{\phi_j^2(T)}{\zeta_j(T)} d T \Big] &= b^{\phi}_j,  \\
\sum_{l=1}^n \Big[ r_l \int_{T_l}^{\hat T_l} \frac{\psi_j^2(T)}{\zeta_j(T)} d T + \hat r_l \int_{\hat T_l}^{T_{l+1}} \frac{\psi_j^2(T)}{\zeta_j(T)} d T \Big] &= b^{\psi}_j,  \qquad j = 1, \ldots, n,
\end{align*}
with right hand side
\begin{align*}
b_j^{\phi} &= \int_0^{T_L} \Big[\mathfrak{r}^{(n)}(T) \frac{\phi_j^2(T)}{\zeta_j(T)} + 
\hat{\mathfrak{r}}^{(n)}(T) \zeta(T) \hat \phi^2_j(T) \Big] dT, \\
b_j^{\psi} &= \int_0^{T_L} \Big[\mathfrak{r}^{(n)}(T) \frac{\psi_j^2(T)}{\zeta_j(T)} + 
\hat{\mathfrak{r}}^{(n)}(T) \zeta(T) \hat \psi^2_j(T) \Big] dT, \qquad j = 1, \ldots, n.
\end{align*}
\end{itemize}
\item \textbf{Output:} The estimates $\zeta^{(n)}(T)$ and $r^{(n)}(T)$.
\end{itemize}
\end{alg}

\subsubsection{A simple estimate of the loss}

We explain in  Appendix \ref{ap:D} that the operator $\cL_\zeta$  defined in \eqref{eq:pencilP} 
is related via a similarity transformation to   the first order Schr\"{o}dinger operator
$\cL + \cQ(T)$ defined in \eqref{eq:F21} and \eqref{eq:SCF3}. In particular, 
\begin{equation}
\begin{pmatrix} 
\zeta^{-\frac{1}{2}}(T) & 0 \\ 0 & \zeta^{\frac{1}{2}}(T) 
\end{pmatrix} \boldsymbol{\Phi}_j^\pm(T) = \frac{1}{\sqrt{2}} \begin{pmatrix} \frac{\phi_j(T)}{\sqrt{\zeta(T)}} \\ \\
\pm i \sqrt{\zeta(T)} \hat \phi_j(T) \end{pmatrix},
\label{eq:EigLQ}
\end{equation}
and similar for $\boldsymbol{\Psi}_j^\pm(T)$ are the eigenfunctions 
$\cL + \cQ(T)$. In the case of a constant $\cQ(T)$ i.e., 
$\frac{d}{dT} \ln \zeta(T) = $ constant, these eigenfunctions have the explicit expression
\begin{align}
\frac{1}{\sqrt{2}} \begin{pmatrix} \frac{\phi_j(T)}{\sqrt{\zeta(T)}} \\ \\
\pm i \sqrt{\zeta(T)} \hat \phi_j(T) \end{pmatrix} = \frac{1}{\sqrt{T_L}} \begin{pmatrix} 
 \cos\big[\frac{(j-1/2) \pi T}{T_L} \big] \\\\
 \sin\big[\frac{(j-1/2) \pi T}{T_L} \big] \end{pmatrix} \label{eq:EIGC1}
 \end{align}
 and 
 \begin{align}
\frac{1}{\sqrt{2}} \begin{pmatrix} \frac{\psi_j(T)}{\sqrt{\zeta(T)}} \\ \\
\pm i \sqrt{\zeta(T)} \hat \psi_j(T) \end{pmatrix} = \frac{1}{\sqrt{T_L}} \begin{pmatrix} 
 \sin\big[\frac{ j\pi T}{T_L} \big] \\\\
 \cos\big[\frac{j \pi T}{T_L} \big] \end{pmatrix}. \label{eq:EIGC2}
 \end{align}
Substituting these expressions in Proposition \ref{prop.3} and using trigonometric identities we get that in this case
\begin{align*}
\int_0^{T_L} \big[ r(T) - \mathfrak{r}^{(n)}(T) - \hat{\mathfrak{r}}^{(n)}(T) \big] d T + 
\int_0^{T_L} \big[ r(T) - \mathfrak{r}^{(n)}(T) + \hat{\mathfrak{r}}^{(n)}(T) \big] \nonumber \\
\times  \cos \Big[\frac{(2j-1) \pi T}{T_L} \Big] dT \approx 0, \qquad  \\
\int_0^{T_L} \big[ r(T) - \mathfrak{r}^{(n)}(T) - \hat{\mathfrak{r}}^{(n)}(T) \big] d T - 
\int_0^{T_L} \big[ r(T) - \mathfrak{r}^{(n)}(T) + \hat{\mathfrak{r}}^{(n)}(T) \big] \nonumber \\
\times  \cos \Big[\frac{2j \pi T}{T_L} \Big] dT \approx 0, \qquad j \ge 1.
\end{align*}
We conclude that 
\begin{equation}
r_0 = \frac{1}{T_L} \int_0^{T_L} r(T) dT \approx \frac{1}{T_L} \int_0^{T_L} \big[\mathfrak{r}^{(n)}(T) + \hat{\mathfrak{r}}^{(n)}(T) \big] d T,
\end{equation}
and using that  $\Big( \cos\big(j \pi T/T_L) \Big)_{j \ge 1}$ is an orthogonal basis in the subspace of mean zero functions in 
$L^2([0,\tau_L])$, we have 
\begin{equation}
r(T) \approx r^{(n)}(T) = \mathfrak{r}^{(n)}(T) - \hat{\mathfrak{r}}^{(n)}(T).
\label{eq:simpleEst}
\end{equation}
This is a simple and explicit inversion formula that does not involve solving a linear system as in Algorithm \ref{alg:invGrid}.

For a more general  impedance function Algorithm 
\ref{alg:invGrid} gives a better result, as illustrated with numerical simulations in the next  section. However, the estimate 
\eqref{eq:simpleEst} is not too far off, because the expressions \eqref{eq:EIGC1}--\eqref{eq:EIGC2} of the eigenfunctions of the Schr\"{o}dinger 
operator $\cL + \cQ(T)$ hold asymptotically, for $j \gg 1$.

\subsection{Optimization approach to inversion and numerical results}
\label{sect:ROMinv3}
In this section we present numerical results and an inversion algorithm based on optimization. The setup for the numerical simulations is described in Appendix \ref{ap:F}. We begin with illustrations of the  inversion on the spectrally matched grid, and then give the optimization method for the more difficult cases that are not addressed in sections \ref{sect:ROMinv1} and \ref{sect:ROMinv2}. We end in section
\ref{sect:ROMinv3noise} with inversion results for noisy data.

\subsubsection{Inversion on the spectrally matched grid}
\begin{figure}[h]
\centering
\includegraphics[clip,trim=0mm 0mm 0mm 0mm,height = 0.3\textwidth, width=0.75\textwidth]{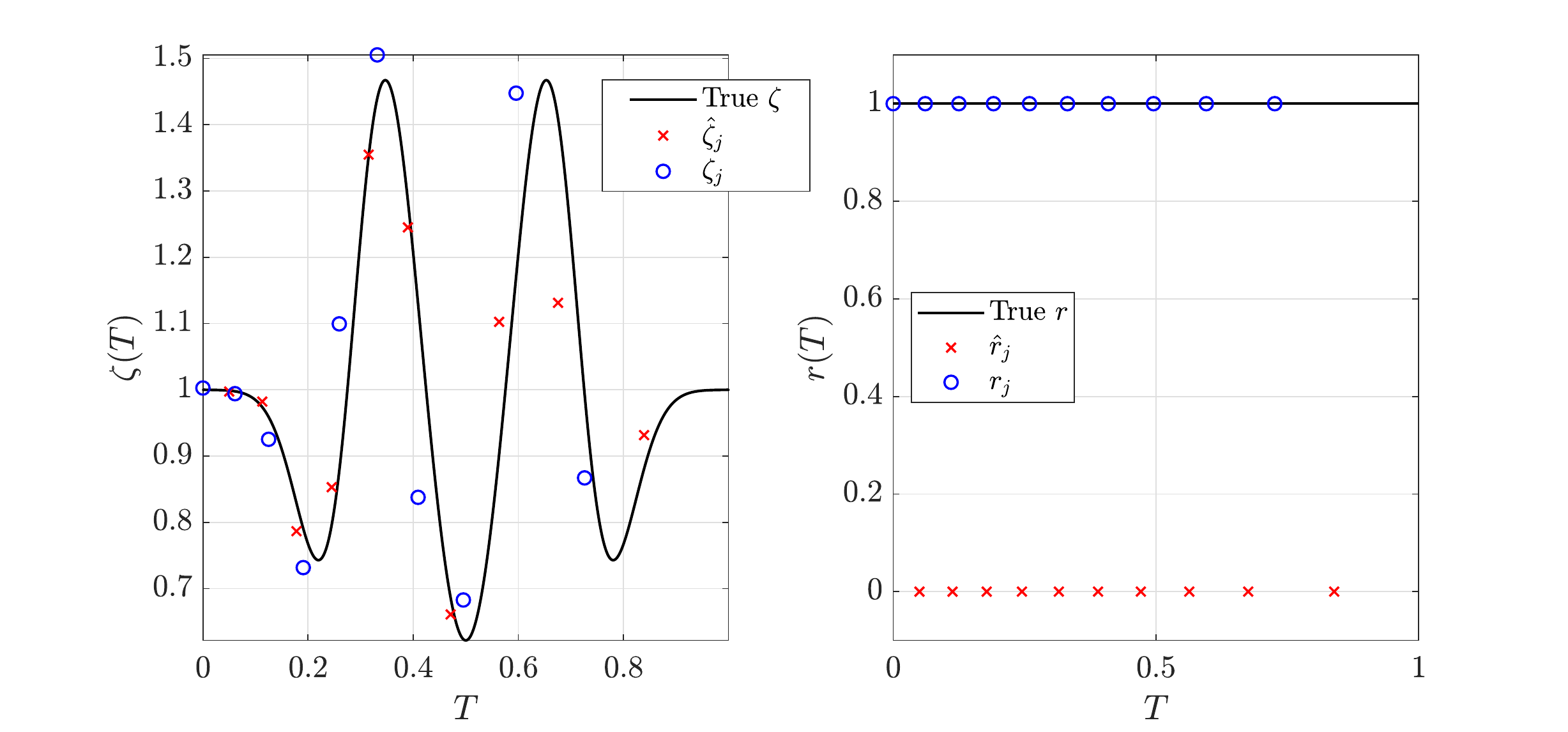} \\
\includegraphics[clip,trim=0mm 0mm 0mm 0mm,height = 0.3\textwidth, width=0.75\textwidth]{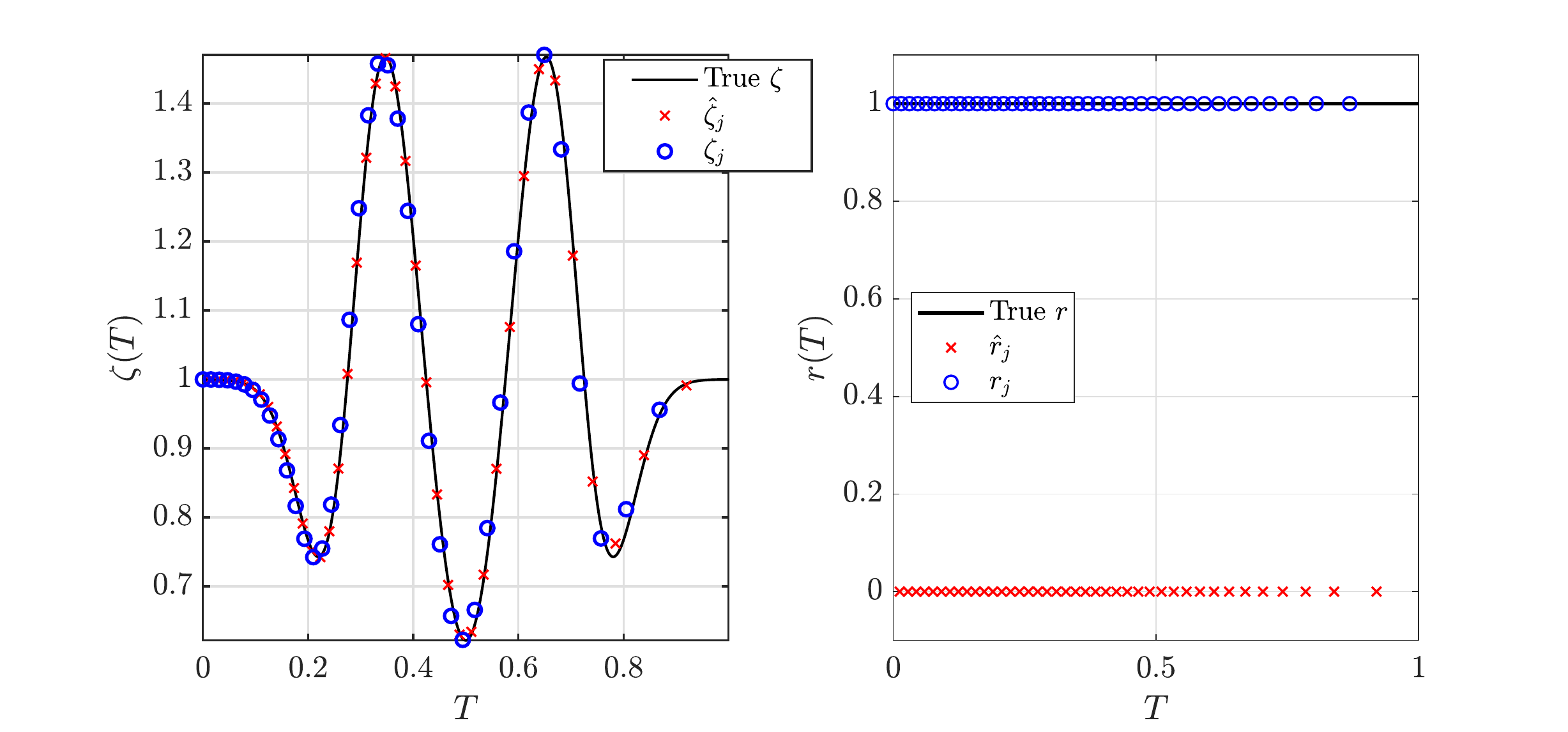}
\vspace{-0.12in}\caption{ROM based inversion on the spectrally matched grid for an impedance function $\zeta(T)$ shown 
in the left plots with the black line. These plots show the values $(\zeta_j)_{j=1}^n$ with the blue circles and 
$(\hat \zeta_j)_{j=1}^n$  with the red crosses.  The right plots show that the constant loss $r_0 = 1$ estimate is exact. 
The values $(\mathfrak{r}_j)_{j=1}^n$ are shown with the blue circles and $(\hat{\mathfrak{r}}_j)_{j=1}^n$ with the red crosses.
The top plots are  for  $n = 10$ and the bottom plots for $n = 40$. }
\label{fig:rec1}
\end{figure}

Propositions \ref{prop.1} and \ref{prop.2} state that in the case of a constant loss $r(T) = r_0$, the ROM coefficients give exactly 
this loss i.e., 
$\mathfrak{r}_j = r_0$ and $\hat{\mathfrak{r}}_j = 0$ for $j = 1, \ldots, n$, and the impedance estimate $\zeta^{(n)}(T)$ defined as in 
Algorithm \ref{alg:invGrid} converges pointwise to $\zeta(T)$ as $n \to \infty$. This is illustrated in Fig. \ref{fig:rec1}.

\begin{figure}[t]
\centering
\includegraphics[clip,trim=0mm 0mm 0mm 0mm,width=\textwidth]{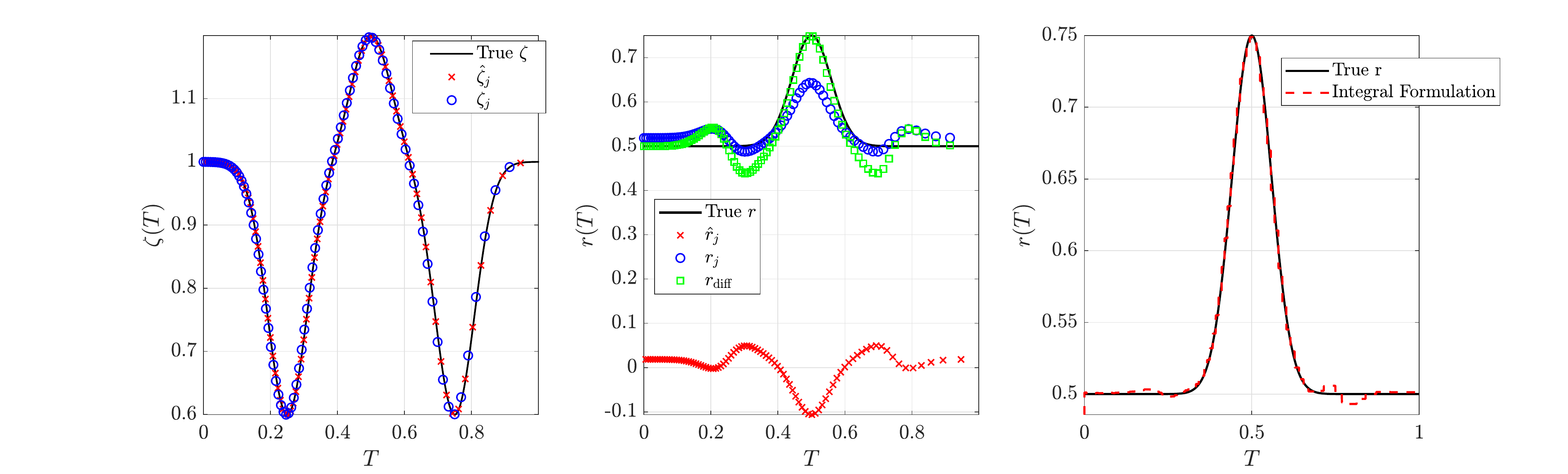}\\
\includegraphics[clip,trim=0mm 0mm 0mm 0mm,width=\textwidth]{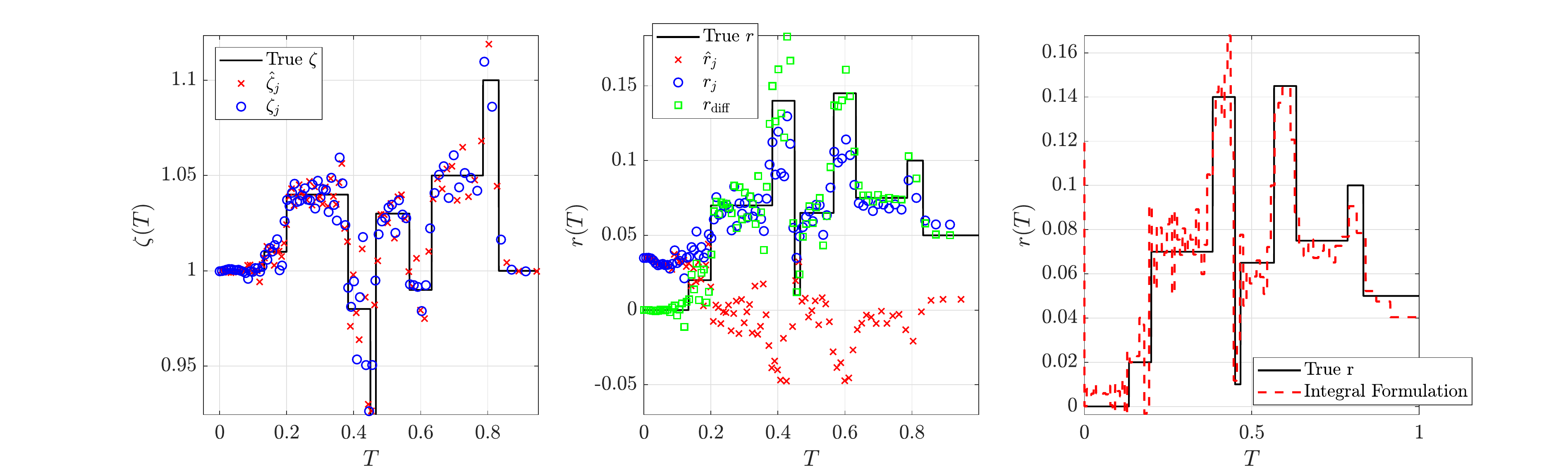}
\vspace{-0.2in}\caption{ROM based inversion on the spectrally matched grid for $n = 90$. The impedance function $\zeta(T)$ 
is shown with the black line 
in the left plot and the loss function $r(T)$ is shown in the middle and right plots with the black line. The left plots show the values $(\zeta_j)_{j=1}^n$  with the blue circles and $(\hat \zeta_j)_{j=1}^n$ 
with the red crosses.  The middle plots show the values $(\mathfrak{r}_j)_{j=1}^n$ 
with the blue circles and $(\hat{\mathfrak{r}}_j)_{j=1}^n$ with the red crosses. The simple estimate \eqref{eq:simpleEst}
is shown with the green circles. The right plots show with the red line the loss estimated with Algorithm \ref{alg:invGrid}. }
\label{fig:rec2}
\end{figure}

In Fig. \ref{fig:rec2} we show the inversion results for  variable loss functions. In the top plots $\zeta(T)$ and $r(T)$ are smooth, as assumed in the analysis. In the bottom plots they are discontinuous.
The impedance function is recovered 
as well as in the previous example (top left plot), but as stated in  Remark~\ref{rem:3}, the ROM coefficients $(\hat{\mathfrak{r}}_j)_{j=1}^n$ are no longer zero.  We display (with green circles) the simple estimate \eqref{eq:simpleEst} of the  loss function and note that it is not too far from $r(T)$. However, the estimate calculated by Algorithm  \ref{alg:invGrid} is much better. The inversion 
for the rougher impedance and loss functions (bottom plots) are not as good as in the smooth medium, due to Gibbs-like  oscillations near the jumps. This is typical of any inversion algorithm and can be mitigated by introducing some regularization in the 
inversion. 

\subsubsection{Optimization approach and numerical results}
\label{sect:numerics}
\begin{figure}[t]
\centering
\includegraphics[width=1.0\textwidth]{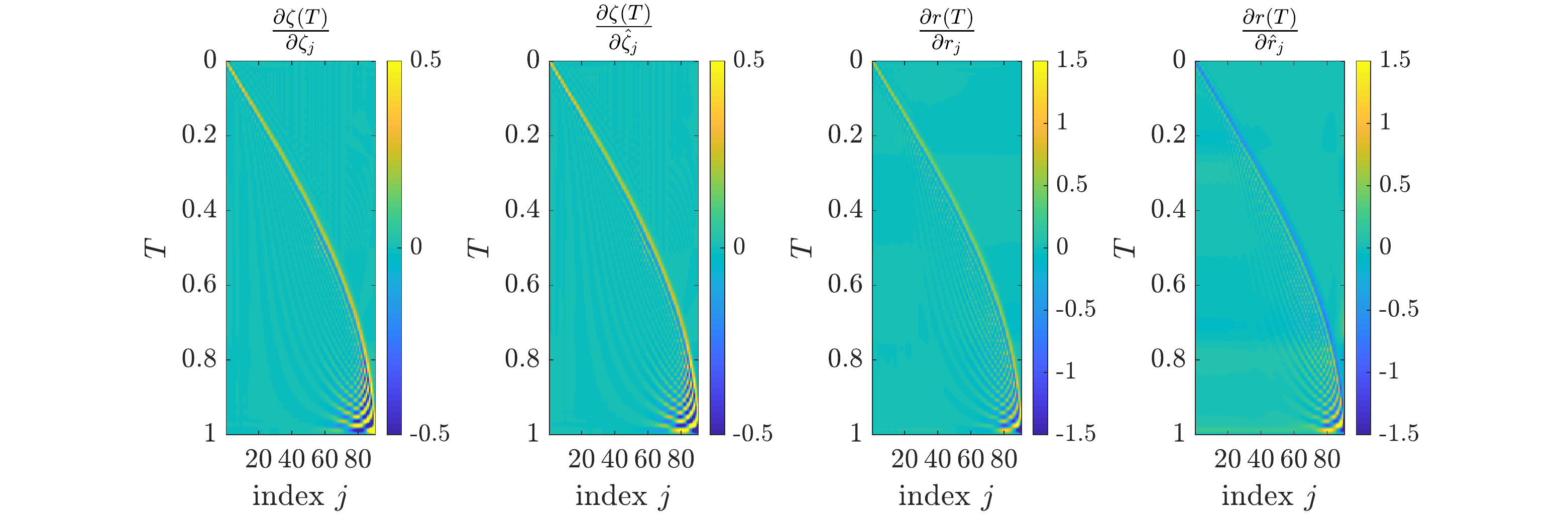}
\vspace{-0.2in}\caption{Sensitivity functions calculated about the true impedance and loss as in Fig. \ref{fig:rec5}.}
\label{fig:rec4}
\end{figure}

If the loss function has larger variations, the estimates given by Algorithm \ref{alg:invGrid} are not so accurate. However,
the ROM coefficients still contain information about $\zeta(T)$ and $r(T)$ that is localized in space, as can be seen  in Fig. \ref{fig:rec4}.  The left two plots in this figure  display the derivative of  the parameters $(\zeta_j, \hat \zeta_j)_{j=1}^n$ with respect to the impedance function $\zeta(T)$. Note how for each index $j = 1, \ldots, n = 90$ they are peaked in magnitude around a specific travel time $T$, which corresponds approximately to the grid points $T_j$ and $\hat T_j$. The right two plots in Fig. \ref{fig:rec4} show a similar behavior of the  sensitivity of the parameters $(\mathfrak{r}_j,
\hat{\mathfrak{r}}_j)_{j=1}^n$ with respect to variations of the loss function $r(T)$.

\begin{figure}[h]
\centering
\includegraphics[clip,trim=0mm 0mm 0mm 0mm,width=\textwidth]{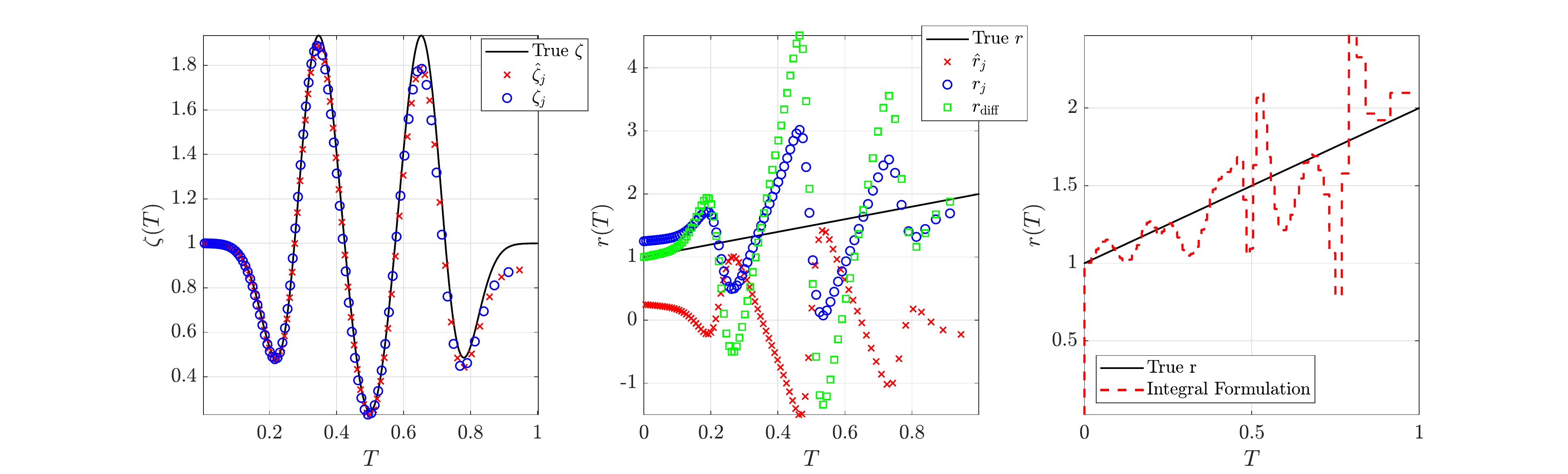}
\includegraphics[clip,trim=0mm 0mm 0mm 0mm,width=0.7\textwidth]{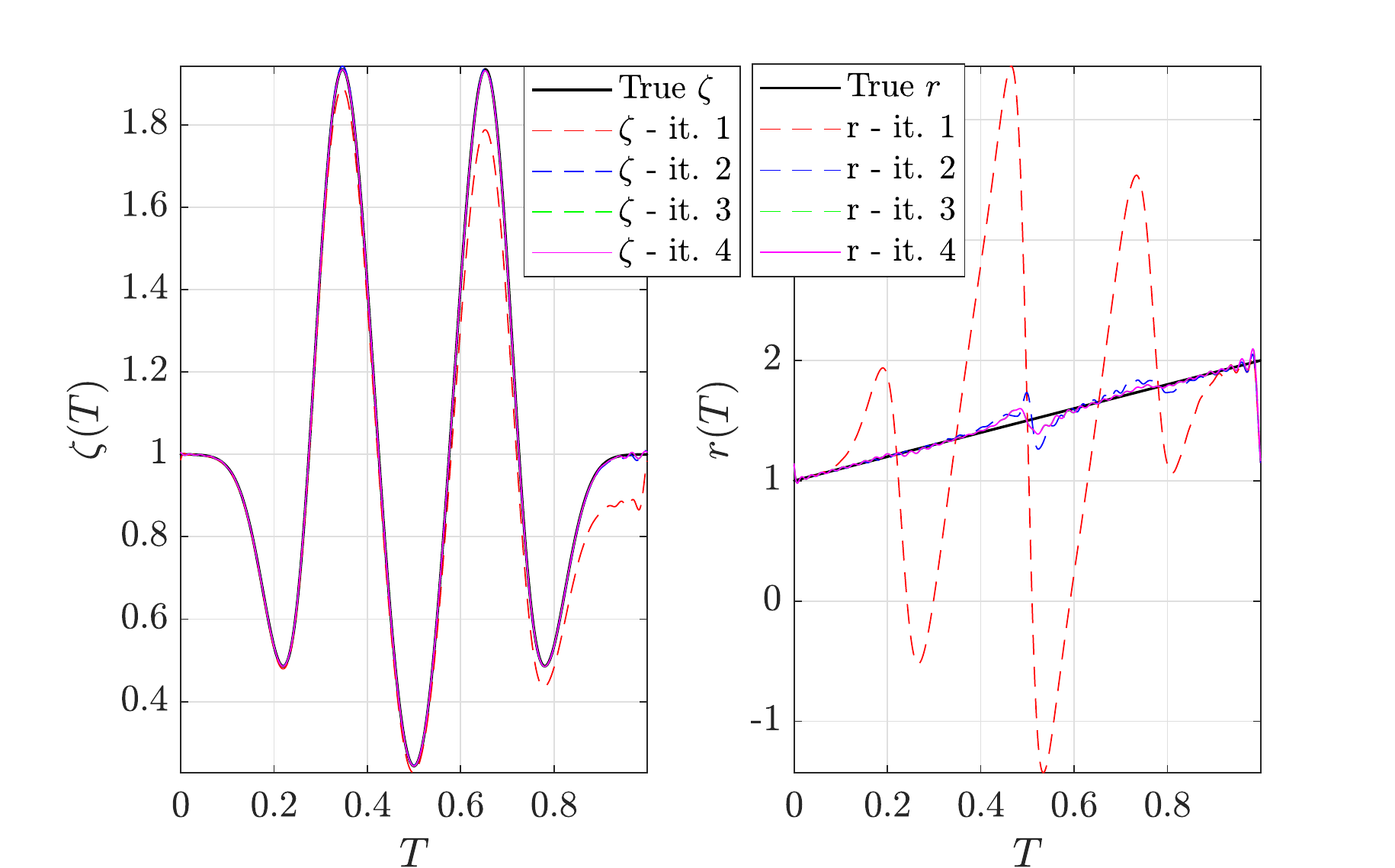}
\vspace{-0.05in}\caption{Top plots show the results obtained with Algorithm \ref{alg:invGrid}. The impedance  is recovered 
well (left plot) but the loss is not. The simple estimate \eqref{eq:simpleEst} is shown with the green circles in the middle plot 
and the estimate calculated by  Algorithm \ref{alg:invGrid} is shown with the red line in the right plot. The optimization results are 
shown in the bottom plots. We have an excellent estimate of both impedance and loss in $4$ iterations. }
\label{fig:rec5}
\end{figure}

This localized information contained in the ROM parameters motivates  a new optimization based approach to inversion. Instead of seeking an approximate inverse of the  nonlinear mapping $\{\zeta(T), r(T)\} \mapsto D(s)$ via least squares data fitting, as is usually done, we propose to invert as follows: First, we map the data to the ROM coefficients: 
\[D(s) \mapsto \{\zeta_j, \hat \zeta_j, \mathfrak{r}_j, \hat{\mathfrak{r}}_j, ~
j = 1, \ldots, n\},\] 
using the non-iterative Algorithms \ref{alg:Lancz1}--\ref{alg:Lancz2} and formulas \eqref{eq:R5}.
Then, we estimate $\zeta(T)$ and $r(T)$ by solving the optimization  
\begin{align}
\min_{\zeta^S, r^S} \sum_{j=1}^n \Big[ 
|\mathfrak{r}_j - \mathfrak{r}_j^{S} |^2 + |\hat{\mathfrak{r}}_j - \hat{\mathfrak{r}}_j^{S} |^2 + 
|\zeta_j - \zeta_j^S|^2 + |\hat \zeta_j - \hat \zeta_j^S|^2 \Big].
\label{eq:OPT}
\end{align}
Here $\zeta^S(T), r^S(T)$ are the search impedance and loss functions in a finite dimensional search space, and $\{ \zeta_j^S, \hat{\zeta}_j^S, \mathfrak{r}_j^{S},  \hat{\mathfrak{r}}_j^{S}, ~~ j = 1, \ldots, n \}$ are the ROM coefficients calculated from the transfer function computed in the medium with impedance $\zeta^S(T)$ and loss $r^S(T)$. In the simulations with 
results displayed in Fig. \ref{fig:rec5} we used the Fourier search space 
\begin{equation*} 
\zeta^S(T), r^S(T) \in \mbox{span}\left\{ \cos\Big[ \pi j \Big( \frac{2T}{T_L}-1\Big)\Big], ~ \sin\Big[ \pi j \Big( \frac{2T}{T_L}-1\Big)\Big], ~
j = 0, \ldots, \Big\lfloor \frac{n}{2} \Big\rfloor \right\},
\end{equation*} 
but other finite dimensional spaces can be used. The minimization \eqref{eq:OPT} is carried out with the Gauss-Newton algorithm,
using the Jacobian calculated as explained in appendix \ref{ap:F}.

As we increase the order of the Fourier series, the linear system solved at each Gauss-Newton step becomes ill-conditioned. In Fig. \ref{fig:rec5} we did not use regularization using a Fourier series with $n=90$ terms, and obtained convergence in 4 iterations. This is generic for all the simulations that we ran. However, if we over-parametrized, i.e., increased $n$ further, then regularization would be needed to stabilize the inversion.

\begin{figure}[h]
\centering

\includegraphics[clip,trim=0mm 0mm 0mm 0mm,width=0.95\textwidth]{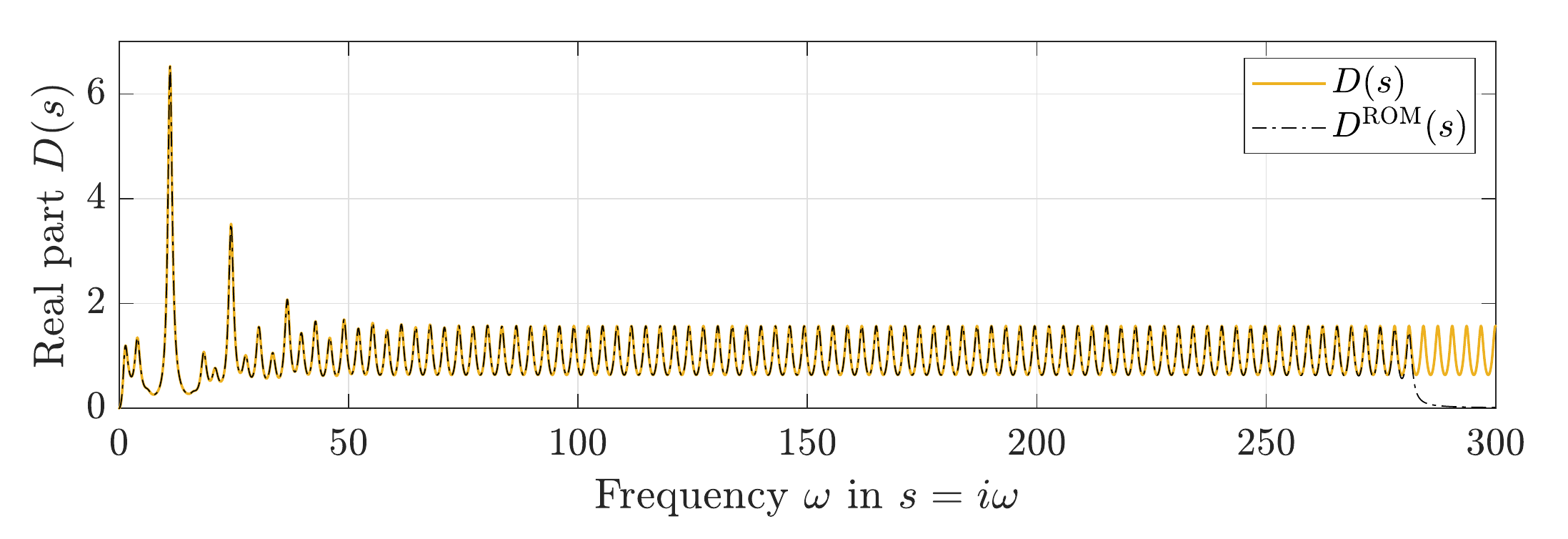}
\vspace{-0.1in}
\caption{Real part of the transfer function  $D(s)$ (full yellow line) and the reduced order model transfer function $D^{\rm ROM}(s)$ (dashed black line) for the configuration shown Figure~\ref{fig:rec5}. The results for the other configurations look very similar. 
That the ROM is passive follows from the fact that $\mbox{Re}\big(D^{\rm ROM}(s)\big) > 0$ for $s$ on the imaginary axis \cite[Corollary 1]{beattie2019robust}.}
\label{fig:rec6}
\end{figure}

\begin{rem}
As showed in appendix \ref{ap:B.p} the dynamical system \eqref{eq:F20} with transfer function \eqref{eq:F23} is passive and 
in pH form. The ROM dynamical system with transfer function $D^{\RM}(s)$ is not in pH form, but it is stable and in all our simulations it is passive. The stability is because the poles of $D^{\RM}(s)$ are the same as the first $n$ poles of $D(s)$, 
and thus lie in the left half complex plane. If $D^{\RM}(s)$ also interpolated the zeroes of $D(s)$, then passivity would be guaranteed by \cite[Theorem 2.1]{sorensen2005passivity}. We only have an approximate interpolation of the zeroes, 
so we can only offer numerical evidence in Figure~\ref{fig:rec6} that  passivity holds.
\end{rem}

\subsection{Inversion for noisy data}
\label{sect:ROMinv3noise}
We repeat here the inversion experiment for the medium with smooth impedance and loss function, shown in the top plots of Fig. \ref{fig:rec2}, using the transfer function contaminated with additive noise denoted by $\mathcal{N}(s)$,
\begin{equation}
D^{\rm noisy}(s) = D(s) + \mathcal{N}(s).
\label{eq:Noise1}
\end{equation}
We refer to appendix Appendix \ref{ap:F} for the details on the sampling of $D(s)$ at $s$ lying in a closed interval on the imaginary axis, centered at the origin. The noise $\mathcal{N}(s)$ is modeled as discrete white Gaussian noise at the sample points in this interval, with 
root mean square equal to $5\%$ of the root mean square of $D(s)$ in the spectral interval of interest.
The truncated spectral measure transfer function $D_n^{\RM}(s)$ of the ROM is estimated from \eqref{eq:Noise1} using the 
rational fitting algorithm ``vectorfit" in \cite{PassiveVecFit,VecFit}. The procedure is described in Appendix \ref{ap:F}, where $D(s)$ should 
be replaced by $D^{\rm noisy}(s)$.

\begin{figure}[t]
\includegraphics[width=\textwidth]{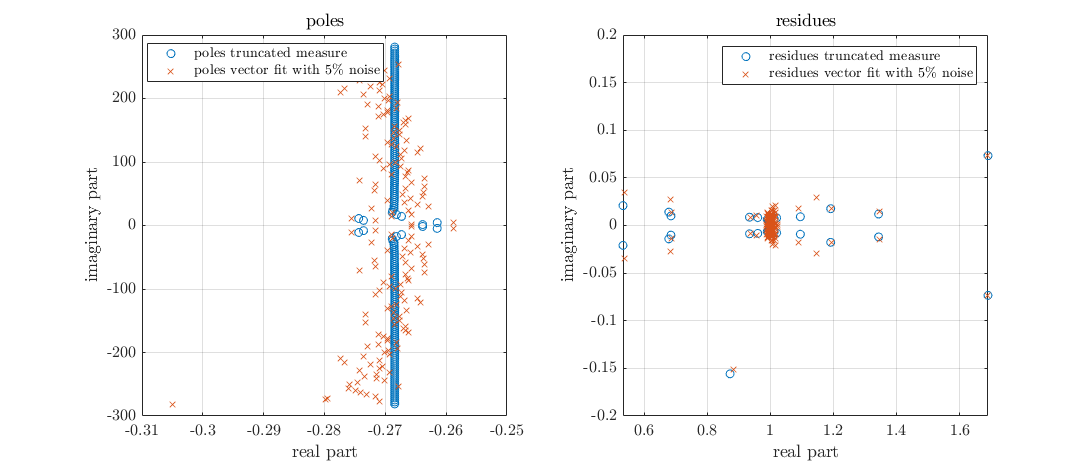}
\caption{Illustration of the poles and residues obtained with vectorfit, for the noiseless and noisy data.}
\label{fig:noise1}
\end{figure}
\begin{figure}[h]
\includegraphics[width=\textwidth]{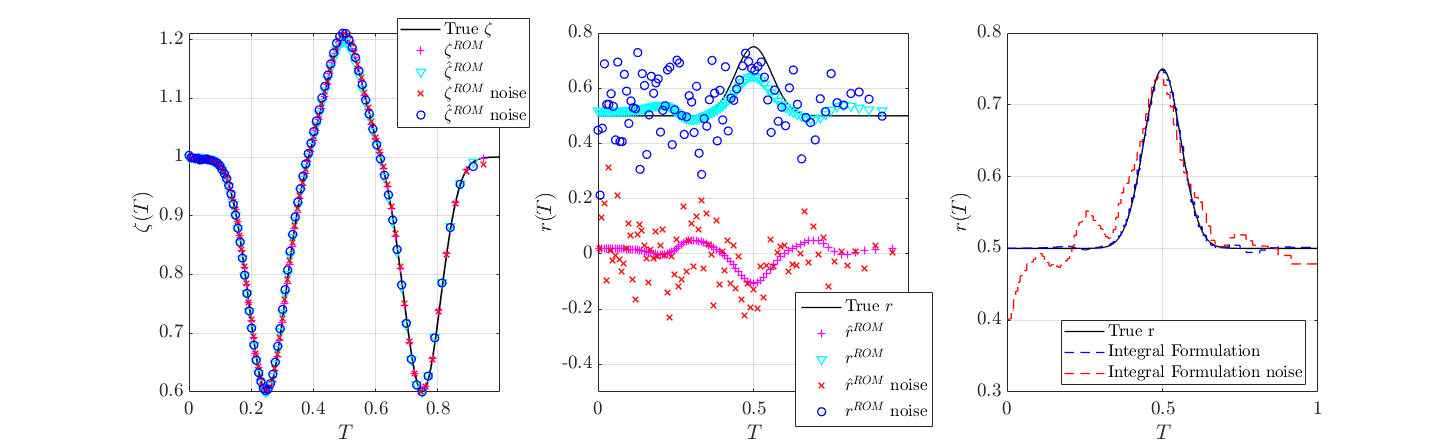}
\caption{The analogues of the top plots in Fig. \ref{fig:rec2} for noisy data. The results with no noise are shown for comparison.}
\label{fig:noise2}
\end{figure}

We illustrate in Fig. \ref{fig:noise1} the poles and residues estimated with the vectorfit algorithm in the spectral 
interval used for inversion. The inversion results are in Fig. \ref{fig:noise2}. Note that to produce the right plot, 
we added regularization penalizing the square of the variation of the loss $r(T)$ to Algorithm \ref{alg:invGrid}.

\section{Summary}
\label{sect:sum}
We introduced a novel approach to inverse scattering in layered lossy media, which is rooted in the reduced order model (ROM) methodology for port-Hamiltonian (pH) dynamical systems. We showed how to transform the system of Maxwell's equations 
for plane waves to a pH dynamical system with transfer function determined by typical remote radar antenna measurements. The dynamical system has two important attributes: it is stable, meaning that the transfer function is analytic in the right half complex plane,  and it is passive, meaning that it does not generate energy internally. 

The ROM is obtained via nonlinear data processing that can be carried out with the computationally efficient J-symmetric Lanczos  algebraic algorithm. It is described by a tridiagonal matrix that looks similar to a three point stencil matrix in a finite difference 
approximation of the dynamical system. It turns out that the entries in this matrix encode local information about the coefficients that define the ROM (the dielectric permittivity $\ep$ and electrical conductivity $\sigma$ of the layered medium), where local means in the cells of a special staggered grid. 
The trouble is that this matrix has non-physical ``magnetic losses" that are not even guaranteed to be positive. Therefore, the ROM dynamical system does not preserve the pH structure, even though it is passive. The main difficulty addressed in the paper is how to embed the ROM in the continuous pH dynamical system, for the 
purpose of solving the inverse scattering problem: Find $\ep$ and $\sigma$ 
from measurements at the antenna, which can be mapped to the transfer function of the pH dynamical system.

The inversion methodology in this paper can be extended to multi-dimensional media. If the measurement setup is as in synthetic aperture radar, where a single antenna emits waves and measures the returns as it moves along some path, one can build a ROM as in this paper for each measurement, and then solve an optimization problem that sums objective functions like \eqref{eq:OPT} over all the ROMs. Such an approach has been used successfully for a parabolic equation in \cite{borcea2014model}.
If the measurements are done by an array of antennas, the transfer function $D(s)$ will be matrix valued. Then, the ROM should be 
obtained via rational approximation of such a transfer function  and it should be brought to block-tridiagonal form via a block Lanczos algorithm. 
This idea was used in \cite{borcea2019reduced} for inverse scattering in lossless media, where the block-tridiagonal 
ROM is used to devise a rapidly converging optimization method. The main difficulty of this approach  is  the absence of scalable algorithms for 
constructing passive, block structured ROMs. Promising recent advances in \cite{benner2020identification,sorensen2005passivity}
have yet to be tried in large scale radar applications.

\begin{acknowledgements}
This research is supported in part by the AFOSR awards FA9550-18-1-0131 and 
 FA9550-20-1-0079, and by ONR  award N00014-17-1-2057. The last two authors of this article are working with Rob Remis and Murthy Guddati on the related project ``Krein's embedding of dissipative data-driven reduced order models''. We would like to thank them for stimulating discussions. We also thank Alex Mamonov and Mikhail Zaslavsky for stimulating discussions on the topic of this paper. \end{acknowledgements}

\vspace{0.1in}
\appendix
\section{From the Weyl function to the transfer function}
\label{ap:A}
The Weyl function $\cW(s)$ is defined\footnote{The Weyl function is denoted by $M$ in \cite{buterin2012inverse} and the Laplace frequency by $\rho$.} in \cite{buterin2012inverse}
as 
\begin{equation}
\cW(s) = \partial_T \phi(0,s),
\label{eq:A2}
\end{equation}
where $\phi(T,s)$ is the solution of 
\begin{equation}
\mathscr{L}_{q,r}(s) \phi(T,s) = 0, ~~\mbox{for}~T \in (0,T_L), \qquad \phi(0,s) = 1, \quad \phi(T_L,s) = 0.
\label{eq:A3}
\end{equation}

Let us introduce the following, pairwise linearly independent solutions associated with the operator pencil 
\eqref{eq:A1}: $\psi(T,s)$, $\xi(T,s)$ and $\eta(T,s)$. The first one satisfies 
\begin{equation}
\mathscr{L}_{q,r}(s) \psi(T,s) = 0, ~~\mbox{for}~T \in (0,T_L), \qquad \psi(T_L,s) = 0, \quad \partial_T \psi(T_L,s) = -1,
\label{eq:A4}
\end{equation}
the second one satisfies 
\begin{equation}
\mathscr{L}_{q,r}(s) \xi(T,s) = 0, ~~\mbox{for}~T \in (0,T_L), \qquad \xi(0,s) = 0, \quad \partial_T \xi(0,s) = 1,
\label{eq:A5}
\end{equation}
and the third one satisfies
\begin{equation}
\mathscr{L}_{q,r}(s) \eta(T,s) = 0, ~~\mbox{for}~T \in (0,T_L), \qquad \eta(0,s) = 1, \quad \partial_T \eta(0,s) = 0.
\label{eq:A6}
\end{equation}
It is easy to check that the Wronskian 
\begin{equation}
\mW_{\psi,\xi}(T,s) = \psi(T,s) \partial_T \xi(T,s) - \xi(T,s) \partial_T \psi(T,s) 
\label{eq:A7}
\end{equation}
is constant in $T$, so we can define 
\begin{equation}
\Delta_D(s) = \mW_{\psi,\xi}(T,s) =\mW_{\psi,\xi}(0,s) = \psi(0,s),
\label{eq:A8}
\end{equation}
where we used the boundary condition in \eqref{eq:A5}. 
Similarly, the Wronskian 
\begin{equation}
\mW_{\eta,\psi}(T,s) = \psi(T,s) \partial_T \eta(T,s) -\eta(T,s) \partial_T \psi(T,s) 
\label{eq:A10}
\end{equation}
is constant in $T$, so we can define 
\begin{equation}
\Delta_N(s) = \mW_{\eta,\psi}(T,s) = \mW_{\eta,\psi}(0,s) = -\partial_T \psi(0,s) ,
\label{eq:A11}
\end{equation}
where we used the boundary condition in \eqref{eq:A6}.

Now it follows that the solution of \eqref{eq:A3} can be written as 
\begin{equation}
\phi(T,s) = \frac{\psi(T,s)}{\Delta_D(s)},
\label{eq:A12}
\end{equation}
and the solution $w(T,s)$ of the Schr\"{o}dinger problem \eqref{eq:S2}--\eqref{eq:S4} is  
\begin{equation}
w(T,s) =  s \zeta_0 \frac{\psi(T,s)}{\Delta_N(s)}, \qquad T \in (0,T_L).
\label{eq:A14}
\end{equation}
The latter is because $\mathscr{L}_{q,r}(s) w(T,s) = 0$ for $T \in (0, T_L)$ by construction, 
and at the boundary we have 
\begin{equation}
\partial_T w(0,s) = - s \zeta_0, \qquad w(T_L,s) = 0.
\label{eq:AAB}
\end{equation}
Now using  \eqref{eq:S4} we obtain the jump condition 
\begin{equation}
\partial_T w(0,s) - \partial_T w(0-,s) = - s \zeta_0,
\end{equation}
which corresponds to the Dirac delta forcing $-s \zeta_0 \delta(T)$ in \eqref{eq:S2}.

Solving for $\psi(T,s)$ in  \eqref{eq:A12} we get 
\begin{equation}
w(T,s) = s \zeta_0 \frac{\Delta_D(s)}{\Delta_N(s)} \phi(T,s),
\label{eq:A15}
\end{equation}
and since $\phi(0,s) = 1$, the transfer function has the expression
\begin{equation}
D(s) = w(0,s) =  s \zeta_0 \frac{\Delta_D(s)}{\Delta_N(s)}.
\label{eq:A16}
\end{equation}
Moreover, taking the $T$ derivative in \eqref{eq:A15}  at $T = 0$ and using the definition \eqref{eq:A2} of the Weyl function and the boundary condition \eqref{eq:AAB}, we obtain that 
\begin{equation}
\cW(s) = -\frac{\Delta_N(s)}{\Delta_D(s)} .
\label{eq:A17}
\end{equation}
This proves equation \eqref{eq:S7}.

The poles of the transfer function are the zeroes of the Weyl function and therefore of the Wronskian \eqref{eq:A11}. They correspond to the eigenvalues of the quadratic operator pencil $\mathscr{L}_{q,r}(s)$ with domain 
$\cS_N$ defined in \eqref{eq:S11}. The zeroes of the transfer function are $s = 0$ and the set of poles of the 
Weyl function, which are the eigenvalues of the quadratic operator pencil $\mathscr{L}_{q,r}(s)$ with domain 
$\cS_D$ defined in \eqref{eq:S10}.

\section{The transfer function for small variations of  the loss function}
\label{ap:B}
To use the perturbation theory in \cite{kato}, consider the first order system formulation \eqref{eq:SCF1}--\eqref{eq:SCF2}. We made the assumption \eqref{eq:S1} to simplify the boundary conditions satisfied by 
$w(T,s)$. Similarly, we shall assume in this appendix that 
\begin{equation}
\frac{d}{d T } \zeta(T_L) = 0,
\label{eq:B2}
\end{equation}
which implies that $\hat w(T,s)$ satisfies a homogeneous Neumann boundary condition at $T = T_L$. 
The loss function has small variations, so we model it as 
\begin{equation}
r(T) = r_0 + \alpha \rho(T),  \qquad \sup_{T \in (0,T_L)}|\rho(T)|/r_0 = O(1),  \qquad 0 < \alpha \ll 1. \label{eq:smallloss}
\end{equation}
Using \eqref{eq:smallloss} in \eqref{eq:SCF1} we obtain 
\begin{equation}
\left[\cL + Q(T) + R_\alpha(T) + s I \right] \begin{pmatrix}  w(T,s) \\ \hat w(T,s) \end{pmatrix} = \begin{pmatrix} 
\zeta_0 \delta(T) \\ 0 \end{pmatrix}, \qquad T \in (0-,T_L) ,
\label{eq:B3}
\end{equation}
where  $\cL$ is the differential operator \eqref{eq:F21} and $Q(T)$ is the skew-symmetric multiplication operator defined in 
terms of $\zeta(T)$ in \eqref{eq:SCF3}. The loss function \eqref{eq:smallloss} is in the multiplication operator 
\begin{equation}
R_\alpha(T) = R_0 + \alpha \begin{pmatrix} \rho(T) & 0 \\ 0 & 0 \end{pmatrix}, \qquad R_0= \begin{pmatrix} r_0 & 0 \\ 0 & 0 \end{pmatrix}.
\label{eq:B5} 
\end{equation}

Note that $i \big[\cL + Q(T)\big]$ acting on the space of functions satisfying the boundary conditions 
\eqref{eq:SCF2} is a self-adjoint (with respect to the Euclidian inner product) indefinite  differential operator on a bounded interval and thus has a countable set of real valued eigenvalues with  no finite accumulation point \cite[Chapter 7]{coddington1955theory}. The
squares of these eigenvalues $(-\theta_j^2)_{j \ge 1}$ are  the same as the eigenvalues in the Sturm-Liouville problems 
\begin{equation}
\left[\frac{d^2}{d T ^2} - q(T) \right] \varphi_j(T) = - \theta_j^2 \varphi_j(T), \quad 
\frac{d}{d T }\varphi_j(0-) = \varphi_j(T_L) = 0, \label{eq:B6}
\end{equation}
and 
\begin{equation}
\left[\frac{d^2}{d T ^2} - \hat q(T) \right] \hat \varphi_j(T) = - \theta_j^2 \hat \varphi_j(T), \quad 
\hat \varphi_j(0-) = \frac{d}{d T} \hat \varphi_j(T_L) = 0, \label{eq:B7}
\end{equation}
where 
\begin{align}
q(T) &= \left[\frac{d}{d T} \ln \zeta^{-\frac{1}{2}}(T)\right]^2 +\frac{d^2}{d T ^2} \ln \zeta^{-\frac{1}{2}}(T) = 
\zeta^{\frac{1}{2}}(T) \frac{d^2}{d T ^2} \zeta^{-\frac{1}{2}}(T), \\
\hat q(T) &= \left[\frac{d}{d T} \ln \zeta^{-\frac{1}{2}}(T) \right]^2- \frac{d^2}{d T ^2} \ln \zeta^{-\frac{1}{2}}(T).
\end{align}
The Sturm-Liouville theory gives that the eigenvalues are simple. 
Assuming that the eigenfunctions $\varphi_j(T)$ and $\hat \varphi_j$ in \eqref{eq:B6}--\eqref{eq:B7} are normalized as 
\begin{equation}
\int_0^{T_L} \varphi_j(T) \varphi_l(T) d T = \int_0^{T_L} \hat \varphi_j(T) \hat \varphi_l(T) d T  = \delta_{jl},
\label{eq:B13}
\end{equation}
then the eigenfunctions of the operator $i \big[\cL + Q(T)\big]$ are 
\begin{equation}
\boldsymbol{\varphi}_j^{\pm}(T) = 
\frac{1}{\sqrt{2}} \begin{pmatrix} 
\varphi_j(T) \\ \pm i \hat \varphi_j(T) \end{pmatrix},
\label{eq:varphiEig}
\end{equation}
and the eigenvalues are $\pm \theta_j$. 
Equivalently, the eigenvalues of $\cL + Q(T)$ are $\pm i \theta_j$. 

Now, if we consider the operator $\cL + Q(T) + R_0$,  for the  constant loss  $r_0$,  the eigenfunctions 
$\boldsymbol{\varphi}_j^\pm(T)$ are still orthonormal and are determined by the components of \eqref{eq:varphiEig} as follows
\begin{equation}
\boldsymbol{\varphi}_j^+(T) = \frac{1}{\sqrt{2 + r_0/\la_j}} 
 \begin{pmatrix} 
\varphi_j(T) \\ \\i \sqrt{1+r_0/\la_j} \hat \varphi_j(T) \end{pmatrix}, \quad \boldsymbol{\varphi}_j^-(T) = 
\overline{\boldsymbol{\varphi}_j^+(T)}, \qquad j \ge 1,
\label{eq:varphiEigr0}
\end{equation}
where $s = \la_j$ is the root of 
\begin{equation}
s(s+ r_0) = -\theta_j^2. \label{eq:B8}
\end{equation}
Indeed, it is easy to check that with $ \boldsymbol{\varphi}_j^+(T) $ given in \eqref{eq:varphiEigr0},
\begin{equation}
\big[\cL + Q(T) + R_0\big]  \boldsymbol{\varphi}_j^+(T) + \la_j \boldsymbol{\varphi}_j^+(T) = 0,
\end{equation}
is equivalent to 
\begin{equation}
\big[\cL + Q(T)\big] \begin{pmatrix} 
\varphi_j(T) \\ i \hat \varphi_j(T) \end{pmatrix} + \sqrt{\la_j(\la_j + r_0)}  \begin{pmatrix} 
\varphi_j(T) \\ i \hat \varphi_j(T) \end{pmatrix} = 0,
\end{equation}
which leads to \eqref{eq:B8}.

 Finally, if we add the perturbation ${R}_\alpha(T)-R_0$, which is clearly a bounded operator with 
$O(\alpha)$ norm,  we can use the analytic perturbation theory in \cite{kato}  to obtain that the 
eigenvalues and eigenprojections are analytic for $\alpha$ in some vicinity of $0$. 
Thus, we can use these eigenprojections to express the wave $w(T,s)$ as a series 
and obtain the expression \eqref{eq:TRANSF} of the transfer function.

\section{Passivity and pH structure}
\label{ap:B.p}
A dynamical system is called passive if it does not generate energy internally \cite{sorensen2005passivity,willems2007dissipative}.
As stated in  \cite[Section 2]{sorensen2005passivity},  this property is realized if the transfer function satisfies
\begin{enumerate}
\item $D(\overline{s}) = \overline{D(s)}$, for all $s \in \CC$. 
\item $D(s)$ is analytic for $\mbox{Re}(s) > 0$.
\item $D(s) + \overline{D(s)} \ge 0$ for $\mbox{Re}(s) > 0$.
\end{enumerate}
It is obvious from the expression \eqref{eq:F23} of the transfer function that it satisfies the first condition. 
The second condition says that the dynamical system is stable i.e., the poles of $D(s)$ are in the left half complex plane. 
We know from section \ref{sect:TF} that these poles are  $\{\la_j, \overline{\la_j}, ~j \ge 1 \}$, where $-\la_j$ are the 
eigenvalues of the operator $\cL + \cQ(T) + R(T)$ acting on the space of functions satisfying the boundary conditions 
 \eqref{eq:SCF2}. If we let $\bV_j(T)$ be the eigenfunctions, then we have 
 \begin{equation}
 \left[\cL + \cQ(T) + R(T) + \la_j I \right] \bV_j(T) = 0.
 \end{equation}
Taking the real part of the inner product with $\bV_j(T)$ and using that $\cL + \cQ(T)$ are skew-symmetric, we get
\begin{align}
0 &= \mbox{Re} \left\{\int_0^{T_L} \overline{\bV_j(T)^T} \left[\cL + \cQ(T) + R(T) + \la_j I \right] \bV_j(T) d T\right\} \nonumber \\
&=\mbox{Re}(\la_j) \|\bV_j\|_{L^2(0,T_L)}^2 + \int_0^{T_L} \overline{\bV_j(T)^T} R(T) \bV_j(T) d T.
 \end{align}
That $\mbox{Re}(\la_j) \le 0$ follows from this equation and the fact that the diagonal multiplication operator  $R(T)$ is positive semidefinite.

It remains to verify the third condition, which has the following physical interpretation: Since $-u(T,s) \vec{\bf e}_{x_2}$ is the electric field and $\hat u(T,s) \vec{\bf e}_{x_1}$ is the magnetic field, 
the Poynting vector at $T = 0$, which determines the power flow, is 
\[
\frac{1}{2} \mbox{Re} \Big\{ \big[-u(0,s) \vec{\bf e}_{x_2}\big] \times \overline{\big[\hat u(0,s) \vec{\bf e}_{x_1}\big]}\Big\} = 
\frac{1}{2} \mbox{Re} \{ u(0,s) \overline{\hat u(0,s)} \} \vec{{\bf e}}_z = \frac{1}{4} \Big[ D(s) + \overline{D(s)} \Big] \vec{{\bf e}}_z.
\]
Thus, the third condition is equivalent to saying that the power flow is into the medium i.e.,  in the positive range direction.

To check this condition, we use the first order system formulation \eqref{eq:SCF1}  to write
\begin{align*}
D(s) &= w(0,s) = \zeta_0 \int_{0-}^T \big(\delta(T), 0\big) \left[ \cL + \cQ(T) + R(T) + s I\right]^{-1} \begin{pmatrix} 
\delta(T) \\ 0 \end{pmatrix} dT  \\
&= \zeta_0 \int_{0-}^T \big(\delta(T), 0\big) \left[ -\cL - \cQ(T) + R(T) + s I\right]^{-1} \begin{pmatrix} 
\delta(T) \\ 0 \end{pmatrix} d T, 
\end{align*}
where in the second line we took the adjoint of the inverse and recalled that $\cL$ and $\cQ(T)$ are skew-symmetric.
Now we can write 
\begin{align}
D(s) + \overline{D(s)} &= \zeta_0 \int_{0-}^T \big(\delta(T), 0\big) \mathcal{P}(s) \begin{pmatrix} 
\delta(T) \\ 0 \end{pmatrix} d T, 
\label{eq:Passive1}
\end{align}
where the operator 
\begin{align*}
&\mathcal{P}(s) = \left[ -\cL - \cQ(T) + R(T) + s I\right]^{-1} + \left[ \cL + \cQ(T) + R(T) + \overline{s} I\right]^{-1} 
\end{align*}
can be factorized as 
\begin{align*}
\mathcal{P}(s) &= \left[ -\cL - \cQ(T) + R(T) + s I\right]^{-1} \Big\{ \left[\cL+ \cQ(T) + R(T)  + \overline{s}I \right] \\
& \quad +  \left[-\cL- \cQ(T) + R(T) + s I\right]
\Big\}   \left[ \cL + \cQ(T) + R(T) + \overline{s} I\right]^{-1} \\
&= 2 \left[ -\cL - \cQ(T) + R(T) + s I\right]^{-1} \left[ R(T) + \mbox{Re}(s) I \right] \left[ \cL + \cQ(T) + R(T) + \overline{s} I\right]^{-1}. 
\end{align*}
Substituting in \eqref{eq:Passive1} and using that 
\[
\left[ \cL + \cQ(T) + R(T) + \overline{s} I\right]^{-1} \begin{pmatrix} 
\delta(T) \\ 0 \end{pmatrix} = \zeta_0^{-1} \begin{pmatrix} w(T,s) \\ \hat w(T,s) \end{pmatrix},
\]
we get that for $ \mbox{Re}(s) > 0$, 
\begin{align} 
D(s) + \overline{D(s)} &= \frac{2}{\zeta_0} \int_{0-}^T \big(\overline{w(T,s)}, \overline{\hat w(T,s)} \big)  \left[ R(T) + \mbox{Re}(s) I \right] \begin{pmatrix} w(T,s) \\ \hat w(T,s) \end{pmatrix} d T \ge 0,
\end{align}
where the inequality is because $R(T)$ is positive semidefinite. 

This proves that our dynamical system is passive. Checking that it has a pH structure is a straightforward verification of \cite[Definition 3]{beattie2019robust}. 

\section{Derivation of the ROM transfer function}
\label{ap:C}
Multiplying equations \eqref{eq:R1} by $\zeta_j$ and \eqref{eq:R2} by $-\hat \zeta_j$, we get the following linear system 
\begin{align}
\left[ 
\bT + s \, {\rm diag}(1,-1,1,-1, \ldots, -1) 
\right]\bU(s) = \frac{{\bf e}_1}{\hat \gamma_1}, 
\label{eq:prefact}
\end{align}
where $\bT$ is the tridiagonal matrix 
\begin{equation}
\bT = \begin{pmatrix} 
r_1 & \frac{1}{\hat \gamma_1} & 0 & 0 & \ldots & 0 & 0 \\
\frac{1}{\gamma_1} & -\hat r_1 & - \frac{1}{\gamma_1} & 0 & \ldots & 0 & 0 \\
0 & - \frac{1}{\hat \gamma_2} & r_2 & \frac{1}{\hat \gamma_2} & \ldots & 0 & 0 \\
&& &\vdots& \\
0 &0 & 0 & 0 & \ldots & \frac{1}{\gamma_n} & - \hat r_n
\end{pmatrix}.
\end{equation}
We can symmetrize this matrix using the diagonal matrix 
$
\bGa = {\rm diag}(\hat \gamma_1, \gamma_1, \hat \gamma_2, \ldots, \gamma_n),
$ so we rewrite \eqref{eq:prefact} as 
\begin{equation}
\Big[ \tilde{\bT} + s  \rm{diag}(1,-1,1,-1, \ldots, -1)\Big] \bGa^{\frac{1}{2}} \bU(s) = \bGa^{\frac{1}{2}} 
\frac{{\bf e}_1}{\hat \gamma_1} = \frac{{\bf e}_1}{\sqrt{\hat \gamma_1}}, 
\end{equation}
where 
 \begin{equation}
\tilde \bT = \bGa^{\frac{1}{2}}\bT \bGa^{-\frac{1}{2}}= \begin{pmatrix} 
r_1 & \frac{1}{\sqrt{\gamma_1 \hat \gamma_1}} & 0 & 0 & \ldots & 0 & 0 \\
\frac{1}{\sqrt{\gamma_1 \hat \gamma_1}} & -\hat r_1 & - \frac{1}{\sqrt{\gamma_1\hat \gamma_2}} & 0 & \ldots & 0 & 0 \\
0 & - \frac{1}{\sqrt{\gamma_1 \hat \gamma_2}} & r_2 & \frac{1}{\sqrt{\gamma_2 \hat \gamma_2}} & \ldots & 0 & 0 \\
&& &\vdots& \\
0 &0 & 0 & 0 & \ldots & \frac{1}{\sqrt{\gamma_n \hat \gamma_n}} & - \hat r_n
\end{pmatrix}.
\end{equation}
Finally, we can factor out the square root of the diagonal matrix in \eqref{eq:prefact} to get 
\begin{equation}
\big[{\bf A} + s {\bf I}\big] {\rm diag}(1,i, 1, i , \ldots, i) \bGa^{\frac{1}{2}} \bU(s) = \frac{{\bf e}_1}{\sqrt{\hat \gamma_1}},
\label{eq:newEq}
\end{equation}
where 
\begin{equation}
{\bf A} =  {\rm diag}(1,-i, 1, -i , \ldots,-i) \tilde \bT {\rm diag}(1,-i, 1, -i , \ldots,-i) 
\end{equation}
is the matrix given in \eqref{eq:defA}--\eqref{eq:R5}.

The ROM transfer function is 
\begin{align} 
D^{\RM}_n(s) &= {\bf e}_1^T \bU(s) = {\bf e}_1^T \bGa^{-\frac{1}{2}} {\rm diag}(1,-i, 1, -i , \ldots,-i) \big[{\bf A} + s {\bf I}\big]^{-1} \frac{{\bf e}_1}{\sqrt{\hat \gamma_1}} \nonumber \\
&= 
{\bf e}_1^T \big[{\bf A} + s {\bf I}\big]^{-1} \frac{{\bf e}_1}{{\hat \gamma_1}}
\end{align}
as stated in equation \eqref{eq:R9}.

\section{Proof of Lemma \ref{lem.1}}
\label{ap:D}
Let us introduce the weighted inner product 
\begin{equation*}
\lb \boldsymbol{\Phi}, \boldsymbol{\Psi} \rb_{\zeta^{-1},\zeta} = \int_0^{T_L} \boldsymbol{\Phi}^\star(T) \begin{pmatrix}
\zeta^{-1}(T) & 0 \\ 0 & \zeta(T) \end{pmatrix} \boldsymbol{\Psi} (T) d T, \qquad \forall ~ \boldsymbol{\Phi}, \boldsymbol{\Psi} \in \big(L^2([0,T_L])\big)^2.
\end{equation*}
The linear operator $i \cL_\zeta$ defined in \eqref{eq:pencilP} with boundary conditions \eqref{eq:BC1} is self-adjoint with 
respect to this inner product and thus has a countable set of real eigenvalues with no finite accumulation point 
\cite[Chapter 7]{coddington1955theory}. This implies that  $\cL_\zeta$ has purely imaginary eigenvalues.
In fact, $\cL_\zeta$ is related via a similarity transformation to the operator $\cL + Q(T)$ studied in Appendix \ref{ap:B} 
\begin{equation}
\cL + Q(T) = \begin{pmatrix} \zeta^{-\frac{1}{2}}(T) & 0 \\ 0 &  \zeta^{\frac{1}{2}}(T) \end{pmatrix} 
\cL_\zeta  \begin{pmatrix}   \zeta^{\frac{1}{2}}(T) & 0 \\ 0 & \zeta^{-\frac{1}{2}}(T) \end{pmatrix} 
\label{eq:D1}
\end{equation}
so the eigenvalues are the same $\{\pm i \theta_j, ~~ j \ge 1\}$. The eigenfunctions of $\cL + Q(T)$ are the vector 
valued functions \eqref{eq:varphiEig}, which are orthonormal in the Euclidian inner product. These define the eigenfunctions of $\cL_\zeta$ 
\begin{equation}
\boldsymbol{\Phi}_j = \frac{1}{\sqrt{2}} \begin{pmatrix} \phi_j (T) \\ \hat \phi_j(T) \end{pmatrix} 
= \begin{pmatrix} \zeta^{\frac{1}{2}} & 0 \\ 0 & \zeta^{-\frac{1}{2}}(T) \end{pmatrix}  \frac{1}{\sqrt{2}} \begin{pmatrix} \varphi_j (T) \\ \hat \varphi_j(T) \end{pmatrix},\qquad j \ge 1,
\end{equation}
and the orthonormality relations \eqref{eq:ORT1} follow from \eqref{eq:B13}.

The same discussion applies to the linear operator $\cL_\zeta$ with boundary conditions \eqref{eq:BC2}. $~ \Box$

\section{Proof of Proposition \ref{prop.3}}
\label{ap:E}
Recall from section \ref{sect:TF1} that the poles of the transfer function $D(s)$ are the 
eigenvalues $\{\la_j, \overline{\la_j}, ~ j\ge 1\}$ of the operator pencil $\mathscr{L}_{q,r}(s)$ defined in \eqref{eq:A1} acting on the space $\cS_N$ defined in \eqref{eq:S11}, whereas the zeroes of $D(s)$ are the eigenvalues 
$\{\mu_j, \overline{\mu_j}, ~ j\ge 1\}$ of the operator \eqref{eq:A1} acting on the space $\cS_D$ defined in \eqref{eq:S10}. 
Moreover, as explained in section \ref{sect:TF} (recall equations \eqref{eq:SCF1} and \eqref{eq:S2}),  $\mathscr{L}_{q,r}(s)$ is connected to the first order pencil 
\begin{equation}
\mathscr{P}_{\zeta,r}(s)  = \cL_\zeta + R_\alpha(T) + s I, \label{eq:Op1}
\end{equation}
with  $R_\alpha(T)$ defined in \eqref{eq:B5} and $\cL_\zeta$ defined in \eqref{eq:pencilP}  
as follows: First, we have the similarity transformation
\begin{equation}
\cL + Q(T) + R_\alpha(T) +s I = \begin{pmatrix} \sqrt{\frac{\zeta_0}{\zeta(T)}} & 0 \\ 0 & \sqrt{\zeta_0 \zeta(T)} \end{pmatrix} 
\mathscr{P}_{\zeta,r}(s) \begin{pmatrix} \sqrt{\frac{\zeta(T)}{\zeta_0}} & 0 \\ 0 & \frac{1}{\sqrt{\zeta_0 \zeta(T)}} \end{pmatrix}. 
\label{eq:E1}
\end{equation}
Second, the first order system 
\begin{equation}
\Big[\cL + Q(T) + R_\alpha(T) +s I\Big] \begin{pmatrix} w(T,s) \\ \hat w(T,s) \end{pmatrix} = 
{\bf 0}, 
\end{equation}
can be written as the second order equation 
\begin{equation}
\mathscr{L}_{q,r}(s) w(T,s) = 0,
\end{equation}
with potential $q(T)$ defined in \eqref{eq:S6} and the other way around. This implies that $\mathscr{L}_{q,r}(s)$ with domain $\cS_N$ has the same eigenvalues $\{\la_j, \overline{\la_j}, ~ j\ge 1\}$ as $\mathscr{P}_{\zeta,r}(s) $ with boundary conditions \eqref{eq:BC1}. Similarly, $\mathscr{L}_{q,r}(s)$ with domain $\cS_D$ has the same eigenvalues $\{\mu_j, \overline{\mu_j}, ~ j\ge 1\}$ as $\mathscr{P}_{\zeta,r}(s) $ with boundary conditions \eqref{eq:BC2}.

We know that the ROM based estimated functions $\zeta^{(n)}(T)$, $\mathfrak{r}^{(n)}(T)$ and $\hat{\mathfrak{r}}^{(n)}(T)$ define an 
operator pencil 
\begin{equation}
\mathscr{P}_{\zeta^{(n)},\mathfrak{r}^{(n)},\hat{\mathfrak{r}}^{(n)}}(s) = \cL_{\zeta^{(n)}} +  \begin{pmatrix} \mathfrak{r}^{(n)}(T) & 0 \\ 0 & \hat{\mathfrak{r}}^{(n)}(T) \end{pmatrix} + s I \label{eq:Op}
\end{equation}
with the following properties: 

\vspace{0.05in}
\begin{enumerate}
\item The pencils \eqref{eq:Op} and \eqref{eq:Op1} with boundary conditions~(\ref{eq:BC1}) have  the same first $n$ eigenvalues $\la_j$, for $j = 1, \ldots, n$. Moreover, the residues $y_j$, which equal the jumps of the spectral measures,  are also the same, 
for $j = 1, \ldots, n.$ 
\item In the case of boundary conditions (\ref{eq:BC2}), the eigenvalues of \eqref{eq:Op} are approximately equal to 
the zeroes $\mu_j$ of the transfer function. We denote the error in their approximation by $o(1)$ in the limit $n \to \infty$.
\end{enumerate} 

\vspace{0.05in}
Now, using the assumption  \eqref{eq:modelr} on the loss function we can write 
\begin{equation}
\mathscr{P}_{\zeta,r}(s)  = \mathscr{P}_{\zeta,r_0}(s)  + \alpha \begin{pmatrix} \rho(T) & 0 \\ 0 & 0 \end{pmatrix}, \label{eq:Op2}
\end{equation}
which is an $O(\alpha)$ perturbation of $\mathscr{P}_{\zeta,r_0}(s)$.  Since the poles and residues of \eqref{eq:Op} and \eqref{eq:Op2} match to all orders of $\alpha$, we can use  Proposition \ref{prop.2} to obtain from the $O(1)$ matching 
that the $O(1)$ primary loss must be $r_0$ and the dual loss is zero. Moreover, the impedance $\zeta^{(n)}(T)$ approximates 
$\zeta(T)$ as $n \to \infty$. 
Therefore, we can write pointwise in $(0,T_L)$ that 
\begin{align}
\zeta^{(n)}(T) &= \zeta(T)\big[1 + o(1) + O(\alpha)\big], \label{eq:E4} \\
\mathfrak{r}^{(n)}(T) &= r_0 + \alpha \rho^{(n)}(T)\big[1 + o(1) + O(\alpha)\big], \label{eq:E5} \\
\hat{\mathfrak{r}}^{(n)}(T) &=\alpha \hat  \rho^{(n)}(T)\big[1 + o(1) + O(\alpha)\big], \label{eq:E6} 
\end{align}
with functions $\rho^{(n)}(T)$ and $\hat \rho^{(n)}(T)$ independent of $\alpha$.

Because the same constant $r_0$ appears in both the pencils \eqref{eq:Op} and \eqref{eq:Op2}, we can subtract it from both 
problems. This results in a transformation of the eigenvalues, but since these eigenvalues match, they 
will be transformed the same way by the subtraction of $r_0$ i.e., they will still match.  The new pencils are 
\begin{equation}
\cL_{\zeta} + \alpha \begin{pmatrix} \rho(T) & 0 \\ 0 & 0 \end{pmatrix} \quad \mbox{and} \quad 
\cL_{\zeta^{(n)}} + \alpha \begin{pmatrix} \rho^{(n)}(T) & 0 \\ 0 & \hat \rho^{(n)}(T)\end{pmatrix},
\label{eq:E10}
\end{equation}
and their eigenvalues are 
\begin{equation} 
\la_j = i \theta_j + \alpha \delta \la_{j} + O(\alpha^2), \qquad \mu_j =  i \vartheta_j + \alpha \delta \mu_{j} + O(\alpha^2),
 \qquad j \ge 1, \label{eq:E11}
\end{equation}
where $i \theta_j$ and $i \vartheta_j$ are the purely imaginary eigenvalues of the lossless problem (see Lemma \ref{lem.1}).
The $O(\alpha)$ perturbations of these eigenvalues are 
\begin{equation}
\delta \la_{j} = \lb \boldsymbol{\Phi}^+_j, \begin{pmatrix} \rho(T) & 0 \\ 0 & 0 \end{pmatrix} \boldsymbol{\Phi}^+_j \rb_{\zeta^{-1},\zeta} =  \lb \boldsymbol{\Phi}^+_j, \begin{pmatrix} \rho^{(n)}(T) & 0 \\ 0 & \hat \rho^{(n)}(T) \end{pmatrix} \boldsymbol{\Phi}^+_j \rb_{\zeta^{-1},\zeta}, \label{eq:E12}
\end{equation}
and similarly 
\begin{equation}
\delta \mu_{j} = \lb \boldsymbol{\Psi}^+_j, \begin{pmatrix} \rho(T) & 0 \\ 0 & 0 \end{pmatrix} \boldsymbol{\Psi}^+_j \rb_{\zeta^{-1},\zeta} =  \lb \boldsymbol{\Psi}^+_j, \begin{pmatrix} \rho^{(n)}(T) & 0 \\ 0 & \hat \rho^{(n)}(T) \end{pmatrix} \boldsymbol{\Psi}^+_j \rb_{\zeta^{-1},\zeta}, \label{eq:E13}
\end{equation}
where $\boldsymbol{\Phi}_j^{+}$ and $\boldsymbol{\Psi}^+_j,$ for $j \ge 1$ are the orthonormal eigenfunctions  in 
Lemma \ref{lem.1}. Substituting their expressions  in these equations gives \eqref{eq:convR1}--\eqref{eq:convR2}.

Finally, note that if $\zeta^{(n)}(T)$ had an 
$O(\alpha)$ error term, then that would reflect in an imaginary perturbation of the eigenvalues, because $\cL_{\zeta}$ 
and $\cL_{\zeta^{(n)}}$ are skew-symmetric operators with respect to the weighted inner product $\lb \cdot, \cdot \rb_{\zeta^{-1},\zeta}$. However, the perturbations \eqref{eq:E12}--\eqref{eq:E12} are real valued, so the error in the impedance must be 
of higher order in $\alpha$. $ ~ \Box$

\section{Setup for the numerical simulations}
\label{ap:F}

In this appendix we describe how we generate the data used in the inversion results 
in Fig. \ref{fig:rec1}--\ref{fig:rec5}.  We  also give details on the calculation of the Jacobian used in the Gauss-Newton iteration for solving the optimization problem \eqref{eq:OPT}. 

To generate the data, we solve the system  (\ref{eq:F4})-(\ref{eq:F5}), with boundary conditions 
\eqref{eq:F6}, \eqref{eq:F17}, using finite differences on a staggered grid with constant step size $\tau$,  except for the 
first dual step size, as illustrated in the following sketch:

\begin{centering}
\begin{tikzpicture}[decorate]
\filldraw[fill=blue!80!white, draw=black,very thick] (2,3) circle (0.2cm);
 \filldraw[fill=blue!80!white, draw=black,very thick] (4,3) circle (0.2cm);
  \filldraw[fill=blue!80!white, draw=black,very thick] (5.9,3) circle (0.2cm);
  \filldraw[fill=blue!80!white, draw=black,very thick] (7.9,3) circle (0.2cm);
  \filldraw[fill=blue!80!white, draw=black,very thick] (10.9,3) circle (0.2cm);
   \filldraw[fill=blue!20!white, draw=black,very thick] (12.9,3) circle (0.2cm);

 \draw[line width=0.4mm] (2.2,3) -- (3.8,3);
 \draw[line width=0.4mm] (4.2,3) -- (5.7,3);
 \draw[line width=0.4mm] (6.1,3) -- (7.7,3);
 \draw[line width=0.4mm, dashed] (8.1,3) -- (10.7,3);
 \draw[line width=0.4mm] (11.1,3) -- (12.7,3);
 
 \node at (3,3.5) {$\tau$};
 \node at (5,3.5) {$\tau$};
 \node at (7,3.5) {$\tau$};

 \filldraw[fill=red!20!white, draw=black,very thick] (2,2) circle (0.2cm);
 \filldraw[fill=red!80!white, draw=black,very thick] (3,2) circle (0.2cm);
 \filldraw[fill=red!80!white, draw=black,very thick] (5,2) circle (0.2cm);
 \filldraw[fill=red!80!white, draw=black,very thick] (7,2) circle (0.2cm);
  \filldraw[fill=red!80!white, draw=black,very thick] (10,2) circle (0.2cm);
 \filldraw[fill=red!80!white, draw=black,very thick] (11.6,2) circle (0.2cm);

 \draw[line width=0.4mm] (2.2,2) -- (2.8,2);
 \draw[line width=0.4mm] (3.2,2) -- (4.8,2);
 \draw[line width=0.4mm] (5.2,2) -- (6.8,2);
 \draw[line width=0.4mm,dashed] (7.2,2) -- (9.8,2);
  \draw[line width=0.4mm] (10.2,2) -- (11.4,2);

 \draw[line width=0.2mm,dashed] (2,1.5) -- (2,3.3);
 \draw[line width=0.2mm,dashed] (12.9,1.8) -- (12.9,3);
  
 \node at (2.5, 1.7) {$\frac{\tau}{2}$};
  \node at (2, 3.5){$T=0$};
 \node at (13, 1.5){$T=T_L$};
  \node at (4, 3.5){$u_2$};
  \node at (6, 3.5){$u_3$};
    \node at (13., 3.5){$u_{N+1}$};
    \node at (3, 1.5){$\hat{u}_2$};
    \node at (5, 1.5){$\hat{u}_3$};
    \node at (11.6, 1.5){$\hat{u}_{N+1}$};
 \node at (4, 1.7) {${\tau}$};
 \node at (6, 1.7) {${\tau}$};
\end{tikzpicture}
\end{centering}

The primary wave $u(T,s)$ is discretized on the primary grid, at the nodes illustrated in blue and the dual wave $\hat u(T,s)$ 
is discretized on the dual grid, at the nodes illustrated in red. We suppress the $s$ dependence of the discretized waves in 
the illustration. The boundary conditions are imposed at the pale colored nodes. The travel time domain is normalized at $T_L=1$ and we use $N=3000$ grid steps, so that the discretization never falls below 30 points per wavelength.
The derivatives are calculated with the standard, two point forward differentiation rule. The truncated measure data can be calculated using the spectral decomposition of the finite differences matrix. However, to better emulate the measurement process, we evaluate $D(s) = u_1(s) $ in the interval $s \in [-i\omega_{\rm max},i\omega_{\rm max}]$, discretized at $10000$ equidistant points, and then use the vectorfit algorithm \cite{PassiveVecFit,VecFit} to extract the poles and residues of $D(s)$. 
The value of $\omega_{\rm max}$ depends on $n$ and it is given in the following table:
\begin{table}[h!]
\begin{center}
\begin{tabular}{ c | c c c}
n& 10 &  40 & 90\\ 
\hline
 $\omega_{\rm max}$& 93 & 124 &281\\  
 Used in  & Figure~\ref{fig:rec1} (top)  &Figure~\ref{fig:rec1} (bottom)& Figure~\ref{fig:rec2} and Figure~\ref{fig:rec5} \\  
\end{tabular}
\end{center}
\label{table:1}
\end{table}

The vectorfit algorithm alone does not give a good estimate of the truncated spectral measure transfer function $D_n^\RM(s)$ 
from $D(s)$, because the poles and residues outside the spectral
interval of interest have a large contribution, especially when the mean loss $r_0$ is large. 
Thus, we proceed as follows: First, we estimate the mean loss $r_0$ by fitting $D(s)$ at points with $|s| \gg 1$ with the transfer function 
\begin{equation}
D^{\rm as}(s;r_0) = \sum_{j=\left\lfloor \frac{T_L \omega_{\rm max}}{\pi} + \frac{1}{2}\right\rfloor}^\infty  \left[ \frac{y_j^{\rm as}}{s - \la_j^{\rm as}} + 
\frac{\overline{y_j^{\rm as}}}{s - \overline{\la_j^{\rm as}}} \right],
\label{eq:VF1}
\end{equation}
with poles and residues given by the asymptotes of the spectrum: 
\begin{equation}
\la_j^{\rm as} = \frac{i(j-1/2)\pi}{T_L} - \frac{r_0}{2}, \qquad y_j^{\rm as} = \frac{\zeta_0}{T_L}\Big[ 1 + \frac{i r_0 T_L}{2(j-1/2)\pi}\Big], \qquad 
j \gg 1.
\end{equation}
Once we estimate $r_0$, we use the vectorfit algorithm to estimate the poles and residues of 
$D(s) - D^{\rm as}(s;r_0)$, for $s \in [-i\omega_{\rm max},i\omega_{\rm max}]$. These then define the truncated 
spectral measure transfer function $D_n^\RM(s)$ of the ROM, according to equation \eqref{eq:SL4}.

The Jacobian used in the Gauss-Newton iteration for solving problem \eqref{eq:OPT} 
is approximated numerically using finite differences. We perturb  the Fourier coefficients of $\zeta^S(T)$ and $r^S(T)$   by $\Delta_{\rm num}=0.01$ and compute the resulting ROM parameters $\mathfrak{r}^{\rm pert}_j, \hat{\mathfrak{r}} ^{\rm pert}_j,\zeta^{\rm pert}_j$ and $\hat\zeta_j^{\rm pert}$. The entries of the Jacobian corresponding to $\zeta_j^S$ are $(\zeta_j^S- \zeta^{\rm pert}_j ) \Delta_{\rm num}^{-1}$, and similarly for the coefficients corresponding to the loss function.


\section*{Data availability}
Data sharing not applicable to this article as no datasets were generated or analysed during the current study. The results of this article are fully reproducible by following the implementation presented in the appendix.
%
%

\bibliographystyle{spmpsci}      
\bibliography{biblio}

\begin{thebibliography}{10}
\providecommand{\url}[1]{{#1}}
\providecommand{\urlprefix}{URL }
\expandafter\ifx\csname urlstyle\endcsname\relax
  \providecommand{\doi}[1]{DOI~\discretionary{}{}{}#1}\else
  \providecommand{\doi}{DOI~\discretionary{}{}{}\begingroup
  \urlstyle{rm}\Url}\fi

\bibitem{beattie2019robust}
Beattie, C., Mehrmann, V., Van~Dooren, P.: Robust port-hamiltonian
  representations of passive systems.
\newblock Automatica \textbf{100}, 182--186 (2019)

\bibitem{benner2020identification}
Benner, P., Goyal, P., Van~Dooren, P.: Identification of port-hamiltonian
  systems from frequency response data.
\newblock Systems \& Control Letters \textbf{143}, 104741 (2020)

\bibitem{CPAM}
Borcea, L., Druskin, V., Knizhnerman, L.: On the continuum limit of a discrete
  inverse spectral problem on optimal finite difference grids.
\newblock Communications on Pure and Applied Mathematics: A Journal Issued by
  the Courant Institute of Mathematical Sciences \textbf{58}(9), 1231--1279
  (2005)

\bibitem{borcea2014model}
Borcea, L., Druskin, V., Mamonov, A., Zaslavsky, M.: A model reduction approach
  to numerical inversion for a parabolic partial differential equation.
\newblock Inverse Problems \textbf{30}(12), 125011 (2014)

\bibitem{borcea2019robust}
Borcea, L., Druskin, V., Mamonov, A., Zaslavsky, M.: Robust nonlinear
  processing of active array data in inverse scattering via truncated reduced
  order models.
\newblock Journal of Computational Physics \textbf{381}, 1--26 (2019)

\bibitem{borcea2019reduced}
Borcea, L., Druskin, V., Mamonov, A., Zaslavsky, M., Zimmerling, J.: Reduced
  order model approach to inverse scattering.
\newblock SIAM Imaging Sciences \textbf{13}(2), 685--723 (2019)

\bibitem{DtB}
Borcea, L., Druskin, V., Mamonov, A.V., Zaslavsky, M.: Untangling the
  nonlinearity in inverse scattering with data-driven reduced order models.
\newblock Inverse Problems \textbf{34}(6), 065008 (2018).
\newblock \doi{10.1088/1361-6420/aabb16}

\bibitem{bruckstein1985differential}
Bruckstein, A.M., Levy, B.C., Kailath, T.: Differential methods in inverse
  scattering.
\newblock SIAM Journal on Applied Mathematics \textbf{45}(2), 312--335 (1985)

\bibitem{buterin2012inverse}
Buterin, S.A., Yurko, V.A.: Inverse problems for second-order differential
  pencils with dirichlet boundary conditions.
\newblock Journal of Inverse and Ill-posed Problems \textbf{20}(5-6), 855--881
  (2012)

\bibitem{chu2005inverse}
Chu, M., Golub, G.: Inverse eigenvalue problems: theory, algorithms, and
  applications.
\newblock Oxford University Press (2005)

\bibitem{coddington1955theory}
Coddington, E., Levinson, N.: Theory of ordinary differentail equations.
\newblock Differential Equations, McGraw-Hill, New York pp. 16--1022 (1955)

\bibitem{druskin2016direct}
Druskin, V., Mamonov, A.V., Thaler, A.E., Zaslavsky, M.: Direct, nonlinear
  inversion algorithm for hyperbolic problems via projection-based model
  reduction.
\newblock SIAM Journal on Imaging Sciences \textbf{9}(2), 684--747 (2016)

\bibitem{freiling2001inverse}
Freiling, G., Yurko, V.A.: Inverse Sturm-Liouville problems and their
  applications.
\newblock NOVA Science Publishers New York (2001)

\bibitem{gugercin2012structure}
Gugercin, S., Polyuga, R., Beattie, C., Van Der~Schaft, A.:
  Structure-preserving tangential interpolation for model reduction of
  port-hamiltonian systems.
\newblock Automatica \textbf{48}(9), 1963--1974 (2012)

\bibitem{VecFit}
{Gustavsen}, B., {Semlyen}, A.: Rational approximation of frequency domain
  responses by vector fitting.
\newblock IEEE Transactions on Power Delivery \textbf{14}(3), 1052--1061
  (1999).
\newblock \doi{10.1109/61.772353}

\bibitem{PassiveVecFit}
{Gustavsen}, B., {Semlyen}, A.: Enforcing passivity for admittance matrices
  approximated by rational functions.
\newblock IEEE Transactions on Power Systems \textbf{16}(1), 97--104 (2001).
\newblock \doi{10.1109/59.910786}

\bibitem{jacob2012linear}
Jacob, B., Zwart, H.: Linear port-Hamiltonian systems on infinite-dimensional
  spaces, vol. 223.
\newblock Springer Science \& Business Media (2012)

\bibitem{jaulent1982inverse}
Jaulent, M.: The inverse scattering problem for lcrg transmission lines.
\newblock Journal of mathematical physics \textbf{23}(12), 2286--2290 (1982)

\bibitem{Joubert}
Joubert, W.: {Lanczos Methods for the Solution of Nonsymmetric Systems of
  Linear Equations}.
\newblock SIAM Journal on Matrix Analysis and Applications \textbf{13}(3),
  926--943 (1992)

\bibitem{kato}
Kato, T.: Perturbation theory for linear operators, vol. 132.
\newblock Springer Science \& Business Media (2013)

\bibitem{lanczos}
Lanczos, C.: An iteration method for the solution of the eigenvalue problem of
  linear differential and integral operators.
\newblock {J. Res. Nat. Bir. Standards} \textbf{45}, 255--282 (1950)

\bibitem{markus2012introduction}
Markus, A.: Introduction to the spectral theory of polynomial operator pencils.
\newblock American Mathematical Soc. (2012)

\bibitem{Marshall1969synthesis}
{Marshall}, T.: {Synthesis of RLC Ladder Networks by Matrix
  Tridiagonalization}.
\newblock {IEEE Transactions on Circuit Theory} \textbf{16}(1), 39--46 (1969)

\bibitem{Morgan2019ReflectionlessFT}
Morgan, M., Groves, W., Boyd, T.: Reflectionless filter topologies supporting
  arbitrary low-pass ladder prototypes.
\newblock IEEE Transactions on Circuits and Systems I: Regular Papers
  \textbf{66}, 594--604 (2019)

\bibitem{pronska2012spectral}
Pronska, N.: Spectral properties of sturm-liouville equations with singular
  energy-dependent potentials.
\newblock arXiv preprint arXiv:1212.6671  (2012)

\bibitem{pronska2013reconstruction}
Pronska, N.: Reconstruction of energy-dependent sturm--liouville equations from
  two spectra.
\newblock Integral Equations and Operator Theory \textbf{76}(3), 403--419
  (2013)

\bibitem{saad1982lanczos}
Saad, Y.: {The Lanczos biorthogonalization algorithm and other oblique
  projection methods for solving large unsymmetric systems}.
\newblock SIAM Journal on Numerical Analysis \textbf{19}(3), 485--506 (1982)

\bibitem{sorensen2005passivity}
Sorensen, D.: Passivity preserving model reduction via interpolation of
  spectral zeros.
\newblock Systems \& Control Letters \textbf{54}(4), 347--360 (2005)

\bibitem{van2004port}
Van Der~Schaft, A.: Port-hamiltonian systems: network modeling and control of
  nonlinear physical systems.
\newblock In: Advanced dynamics and control of structures and machines, pp.
  127--167. Springer (2004)

\bibitem{willems2007dissipative}
Willems, J.: Dissipative dynamical systems.
\newblock European Journal of Control \textbf{13}(2-3), 134--151 (2007)

\bibitem{yagle1989one}
Yagle, A.E.: One-dimensional inverse scattering problems: an asymmetric
  two-component wave system framework.
\newblock Inverse problems \textbf{5}(4), 641 (1989)

\end{thebibliography}

\end{document}